\documentclass[12pt,reqno]{amsart}
\usepackage{amsmath,amsthm, amssymb}

\subjclass{Primary: 37E05 - Secondary: 37D25, 37D50, 60F05}

\theoremstyle{plain}
\newtheorem*{theorem*}{Theorem}
\newtheorem{theorem}{Theorem}
\newtheorem{proposition}{Proposition}[section]
\newtheorem{lemma}[proposition]{Lemma}

\newtheorem*{proposition*}{Proposition}
\newtheorem{corollary}[proposition]{Corollary}

\theoremstyle{definition}

\newtheorem*{definition*}{Definition}
\newtheorem{conjecture}{Conjecture}
\newtheorem*{question*}{Question}

\theoremstyle{remark}

\newtheorem*{remark*}{Remark}

\def\lineclear{\rule{0pt}{0pt}\par\noindent}
 \def\phi{\varphi}
\def\Crit{\ensuremath{\mathcal C } }
\def\Cor{\ensuremath{\mathcal C}}

\begin{document}

\titlepage
\title{Decay of correlations
in one-dimensional dynamics}
\author{Henk Bruin, Stefano Luzzatto and Sebastian van Strien}
\date{\today}

\maketitle

\begin{abstract}
     We consider multimodal $C^3$ interval maps \( f \) satisfying a
     summability  condition on the derivatives \( D_{n} \)
     along the critical
     orbits which implies the existence of an absolutely
     continuous \( f \)-invariant probability measure \( \mu \).
     If \( f \) is non-renormalizable, \( \mu \) is mixing and we show
     that the speed of mixing (decay of correlations) is strongly related to
     the rate of growth of the sequence $(D_n)$
     as \( n\to \infty \).
     We also give sufficient conditions for $\mu$ to satisfy
     the Central Limit Theorem.
     This applies for example to the quadratic Fibonacci map which is shown
     to have subexponential decay of correlations.
\end{abstract}

\bigskip
\bigskip

\begin{center}
{\bf 
D\'ECROISSANCE DES CORRELATIONS POUR \\
DES APPLICATIONS EN DIMENSION $1$
}
\end{center}

\begin{abstract} 
Nous consid\'erons des applications
$C^3$ et multimodales de l'intervalle
ayant des deriv\'es $D_n$ qui satisfont
une condition de sommabilit\'e  le long de l'orbite
critique, ceci entra\^{\i}nant l'existence d'une mesure de probabilit\'es
$\mu$ absolument continue par rapport \`a la mesure de Lebesgues.
Si $f$ n'est pas renormalisable,
$\mu$ est  m\'elangeante et nous
montrons que le d\'ecroissance de corr\'elation est fortement li\'ee
au rapport de croissance de la suite $(D_n)$
lorsque $n \to \infty$. Nous donnons \'egalement une condition suffisante
pour que $\mu$ satisfasse au Th\'eor\`eme de Limite Centrale.
Ceci implique par exemple que l'application quadratique de Fibonacci
poss\`ede un d\'ecroissance de corr\'elation sous-exponentielle.       
\end{abstract}
\newpage

\tableofcontents

\newpage

\section{Introduction}

\subsection{Statement of results}

Let \( f: I \to I \) be a $C^3$ interval or circle map with
a finite critical set \( \Crit \) and no stable or neutral periodic orbit.
All critical points are assumed to have the same finite
{\em critical order} $\ell \in (1, \infty)$. This means that
for $c \in \Crit$,
there exist a diffeomorphism $\phi:\mathbb R \to \mathbb R$
fixing $0$ such that for $x$ close to $c$,
$$
f(x) = \pm |\phi(x-c)|^{\ell} + f(c),
$$
where the $\pm$ may depend on $\mbox{sgn}(x-c)$.
For a critical point $c$, let
$$
D_{n}(c) = |(f^{n})'(f(c))|.
$$
The aim of this paper is to prove the existence of an
absolutely continuous invariant probability measure
under sufficient growth conditions of $D_n(c)$ and
to study its statistical properties (rate of mixing, Central Limit Theorem).
In the proofs we will use distortion estimates valid for maps with 
negative Schwarzian derivative. By a result of Kozlovski \cite{Koz}
(generalized to the multimodal setting by 
van Strien \& Vargas \cite{vSV}),
similar estimates hold if $f$ is $C^3$ and has no stable or neutral
periodic orbit.

\begin{theorem}[Existence of invariant probability measures]
\label{existence}
If \( f \) satisfies 
\begin{equation}
\sum_n  D_n^{-1/(2\ell-1)}(c) < \infty \mbox{ for each } c \in \Crit, 
\tag{$*$}
\end{equation}
then there exists an \( f \)-invariant probability measure \( \mu \) absolutely
continuous with respect to Lebesgue measure ({\em acip}).
\end{theorem}

Moreover, as was shown in general in \cite{Kel},
$\mbox{supp}(\mu)$ is an interval or cycle of intervals.
In the unimodal case, Nowicki \& van Strien \cite{NS91} proved the same
result under the assumption that
$\sum_n D_n^{-1/\ell} < \infty$.
Theorem~\ref{existence} is the first general existence theorem in the
multimodal case, see \cite{BLS}.
Bruin \& van Strien \cite{BvS2} later proved the following generalization:
if $\sum_n D_n^{-1/\ell_{\max}}(c) < \infty$ for all $c \in \Crit$
and $\ell_{\max} = \max \ell(c)$, then $f$ has an acip.
In the holomorphic case,
Przytycki proved that for every $\alpha$-conformal measure $\mu$ on
$J$ (satisfying an additional assumption),
assuming $D_n(c)$ grows exponentially, there exists an invariant
probability measure which is absolutely continuous with
respect to $\mu$, see \cite{Prz}.
Condition ($*$) is equivalent to the following (see Lemma~\ref{equiv}):
\begin{quote}
There exists a sequence
$\{ \gamma_n \}$,
$0 < \gamma_n < \frac12$, such that
$\sum_n \gamma_n < \infty$ and
\begin{equation}
\sum_n [\gamma_n^{\ell-1} D_n(c)]^{-1/\ell} < \infty \text{ for all }
c \in \Crit. \tag{$**$}
\end{equation}
\end{quote}
We prefer to use this version of ($*$), as the $\gamma_i$ play an
important role in the binding method in the proof, and in the formulation
of the other theorems. Let us also abbreviate
\begin{equation}\label{bn}
b_n(c) := \left[ \gamma_n^{\ell-1} D_n(c)\right]^{-1/\ell}.
\end{equation}

The measure $\mu$ need not be unique if $f$ is multimodal
and not Lebesgue ergodic
(unimodal maps with negative Schwarzian derivative
are Lebesgue ergodic, \cite{BL}).
Let $X$ be a closed $f$-invariant set of positive Lebesgue
measure such that $X$ contains no smaller set with these properties.
Misiurewicz \cite{mis} proved that $X$ contains a critical point.
We do not know if $X$ can be a Cantor set
(however, cf. \cite{V}),
but if $X$ supports the measure $\mu$ from Theorem~\ref{existence},
then $X$ has a non-empty interior.
By a result of Ledrappier \cite{Ledrappier}, $\mu$ is mixing
if and only if $f:X \to X$ is not {\em renormalizable}, i.e.
$X$ is not a cycle of intervals permuted by $f$.
In this case it is natural to ask about the
speed of mixing, quantified through the
\textit{correlation function}
\[
\Cor_{n}=
\Cor_{n}(\phi, \psi) =
 \left|\int(\phi \circ f^{n}) \psi d{\mu} -
    \int \phi d\mu\int\psi d\mu\right|,
\]
where \( \phi \) and \( \psi \) are respectively bounded and H\"older
continuous functions on $X$.

Write
\begin{equation}\label{dn}
d_n(c) := \min_{i < n} [\gamma_i/D_i(c)]^{1/\ell} |f^i(c) - \Crit|.
\end{equation}
Obviously, $d_n(c) \leq \gamma_{n-1} b_{n-1}(c) < b_{n-1}(c)$.

\begin{theorem}[Decay of correlations]\label{decay}
Let $f$ satisfy ($*$) and let $\mu$
be an absolutely continuous invariant probability measure
with support $\mbox{supp}(\mu)$.
If $f$ is not renormalizable on $\mbox{supp}(\mu)$, then
$(\mbox{supp}(\mu), \mu, f)$ is mixing with the following rates:
\begin{description}
  \item [Polynomial case]
 If
 $$
 d_n(c) \leq C n^{-\alpha}
 $$
for all $c \in \Crit$, some $\alpha > 1$ and all $n\geq 1$,
then for each
$$
\tilde \alpha < \alpha - 1.
$$
there exist $\tilde C = \tilde C(\phi,\psi) > 0$ such that
\[
\Cor_n \leq \tilde C  n^{-\tilde\alpha }
\text{ for all } n \geq 1.
 \]
 \item [Stretched exponential case]
 If
 \[
 b_n(c) \leq Ce^{-\beta n^{\alpha}}
 \]
for all $c \in \Crit$,
some $C, \beta > 0$, $\alpha \in (0,1)$ and all $n \geq 1$, then
for all $\tilde \alpha \in (0, \alpha)$ there exist
$\tilde C = \tilde C(\phi,\psi), \tilde \beta > 0$ such that
 \[
 \Cor_{n}\leq \tilde C e^{-\tilde\beta n^{\tilde\alpha}}
\text{ for all } n\geq 1.
 \]
 \item [Exponential case]
 If
 $$
 b_n(c) \leq Ce^{-\beta n}
 $$
for all $c \in \Crit$, some $C, \beta > 0$
and $n\geq 1$, then there exist $\tilde C= \tilde C(\phi,\psi), \tilde \beta
> 0$ such
that
 \[
 \Cor_n \leq \tilde C
e^{-\tilde\beta n} \text{ for all } n\geq 1.
 \]
   \end{description}
\end{theorem}

\bigskip
\noindent
Notice that \( d_{n} \) may decay much more rapidly than the terms of
the series in condition \( (**) \). The formulation
in terms of \( d_{n} \) gives us an edge in the polynomial case. 
As an illustration, let us consider the case:
$D_n(c) \geq Cn^{\tau}, \ \tau > 2\ell - 1$ for all
$c \in \Crit$ and $n \geq 1$.
Theorem~\ref{decay} then tells us
that $\Cor_n \leq \tilde C n^{-\tilde \tau}$
for any $\tilde \tau < \frac{\tau-1}{\ell-1} - 1$.
Another use of the $d_n$'s involves the quadratic Fibonacci map, see 
Corollary~\ref{fibo}

\medskip

If \( \mu \) is an \( f \)-invariant probability measure,
we say that the Central Limit Theorem holds if given a H\"older
continuous function \( \phi \) which is not a coboundary (\( \phi\neq
\psi\circ f - \psi \) for any \( \psi \)) there exists \( \sigma>0 \)
such that for every interval \( J\subset \mathbb R \),
\[
\mu\left\{x\in X:
\frac{1}{\sqrt n}\sum_{j=0}^{n-1}\left(\phi(f^{j}(x))-\int\phi d\mu
\right)\in J \right\} \to
\frac{1}{\sigma \sqrt{2\pi} }\int_{J} e^{-t^{2}/ 2\sigma^{2}}dt.
\]
This property is indicative of a certain regularity in the way
Birkhoff averages of H\"older observable approach their expected
asymptotic values.

\begin{theorem}[Central Limit Theorem]\label{CLT}
Let $f$ satisfy ($*$).
If $f$ is not renormalizable and $d_n(c) \leq C n^{-\alpha}, \ \alpha > 2$,
for all $c \in \Crit$ and $n\geq 1$,
then the measure $\mu$ of Theorem~\ref{existence}
satisfies the Central Limit Theorem.
\end{theorem}

The statements about decay of correlations and Central Limit Theorem
in the unimodal exponential case were proved in \cite{KN92,Y92}.
As far as we know the results in all other cases are new.
Most known examples of systems with strictly subexponential
decay of correlations consist of maps which are
uniformly expanding except for the presence of some neutral
fixed point, see for example \cite{LSV,Y98}. The situation here is more
subtle as the cause for the loss of
exponential estimates is not so localized.

In fact in the unimodal case, $D_n \geq C e^{\beta n}$
if and only if $(X,\mu,f)$
has exponential decay of correlations \cite{NS98}.

It is interesting to apply the results to the {\it Fibonacci maps},
i.e. the conjugacy class of unimodal maps characterized by the property
that the sequence of closest return times is exactly the Fibonacci
sequence. Lyubich \& Milnor \cite{LM} proved that
in the quadratic case, the Fibonacci map satisfies
Nowicki's \& van Strien's summability condition.
Here we show

\begin{corollary}\label{fibo} 
Let $f$ be a Fibonacci map with quadratic
critical point. Then one has (faster than)
polynomial decay of correlations
and the central limit theorem holds.
\end{corollary}

\begin{proof}
In fact, the estimates in \cite[Section 5]{LM} show that
condition ($*$) holds for e.g. $\gamma_i = 0.01 \sqrt{1/D_i}$,
so $b_n = 10 D_n^{-1/4}$.
In \cite[Lemma 5.9]{LM}, it is shown that $\sum_n D_n^{-\alpha} < \infty$
for any $\alpha > 0$, which leads to the existence proof of an acip.
Theorem 2 and 3 also hold. Indeed,
if $S_r \approx \gamma^{-r}$ is the $r$-th Fibonacci number
(with  $\gamma = (\sqrt 5 - 1)/2$), then
$|f^{S_r}(c)-c| \approx e^{-\beta' r^2} \approx S_r^{-\beta \log S_r}$
for some $\beta', \beta > 0$.
Let $S_{r-1} < k \leq S_r$ be arbitrary.
Then $d_k \leq S_{r-1}^{-\beta \log S_{r-1}} \leq
(\gamma k)^{-\gamma \beta \log(\gamma k)}$, which
decreases faster than any polynomial, but more slowly
than what we call stretched exponentially.
In particular, the Central Limit Theorem holds.
\end{proof}

\subsection{Techniques and conjectures}

Our approach is to construct an induced Markov map and apply
the result of L.-S. Young \cite{Y} 
which shows that the decay of correlations is
tightly linked to the \emph{tail estimates} of the inducing times.
However, the construction of Markov induced maps (and the corresponding
tower) is quite involved if the map has critical points.
Expanding Markov induced maps have been constructed before,
but only in the unimodal Collet-Eckmann setting 
tail estimates were undertaken.
For our results, we need a new construction, which can be used for
much weaker growth conditions on the orbits of multiple  critical points,
and indeed enables tail estimates of the inducing times.
Apart from its use for estimating decay of correlations, towers were recently
used by Collet \cite{Col99} to describe return time statistics to small 
neighbourhoods. Indeed, combining our results
(namely the tower structure with exponential tail behaviour, 
cf. Subsection~\ref{final_proofs})
with Collet's paper, we can
conclude that for all Collet-Eckmann multimodal maps with constant 
critical order, the quantity $\sup_{ i \leq n} -\log |x-f^i(y)|$ 
satisfies Gumbel's law for $\mu$-a.e. $x$, see \cite[Theorem 1.1]{Col99}
for details. 

Since the growth of derivatives outside a neighbourhood of the critical
set is exponential, one can argue that
the tail is exponential for intervals which spend most of the time
outside such neighbourhoods. Thus we need to concentrate particularly
on intervals which fall inside critical neighbourhoods.
One of the key ideas is to use a shadowing (or
binding) argument
to compare derivative growth for pieces of orbit
to piece of critical orbit that they shadow.
Binding arguments were developed by Jakobson \cite{Jak} and
Benedicks \& Carleson \cite{BC},
under strong growth assumptions ($D_n \geq e^{\sqrt{n}}$ or even
$D_n \geq e^{\lambda n}$) and slow recurrence of the critical point:
$|f^n(c) - c| \geq e^{\alpha n}$ for some small $\alpha$.
This is the so-called basic assumption of \cite{BC}.
Similar conditions were used in several papers
concerned with dynamical and stability properties  of various classes
of one-dimensional maps. We mention \cite{BaV96} in particular where strong
stochastic stability (for random perturbations) was proved, see also
\cite{Thu98} where some similar conditions are
introduced in the context of maps with completely flat critical points.

We dispense with the slow recurrence assumption altogether, and introduce some
new arguments in the construction:
  \begin{itemize}
\item Our definition of binding period (see \eqref{bind}) incorporates the
recurrence pattern of the critical set. As a result, the partition of
the space into intervals of constant induce time is not
fixed in advance, as is the case in \cite{BC}.
\item In order to still count and measure the lengths of partition 
elements, we need intricate combinatorial counting arguments,
which involves assigning {\em itineraries} to the partition elements, 
which indicate
the ``deepness'' of the successive visits to a neighbourhood of the critical
point.
\item Our inducing time consist of three explicit parts: 
the first part is used to recover from the small derivatives near the critical
set (thus achieving expansion); in the second intervals reach 
``large scale'' and the third part is used to reach a 
prefixed interval.
\end{itemize}
In spite of the many differences, we believe that the construction is 
sufficiently robust as to justify
\begin{conjecture}\label{stochstab}
Multimodal Collet-Eckmann maps are strongly stochastically stable.
\end{conjecture}
Tsujii's result on weak stochastic stability \cite{Tsu92} 
indicate in this direction.
Possibly, the Collet-Eckmann condition itself can be
replaced by a much weaker growth condition.

\medskip

Let \( \Lambda \) be a compact (forward)
invariant set for a
smooth map \( f \) and $\mu$ be an \(f\)-invariant ergodic probability
measure.
The measure \( \mu \) is called \emph{hyperbolic} if all
the  Lyapunov exponents corresponding to \( \mu \) are
non-zero (recall that by Oseledec's Theorem, the Lyapunov exponents
associated to a measure are well defined);
it is called a
{\it physical measure} if the set of \( \mu \)-generic points
has positive probability with respect to the given
reference (Lebesgue) measure.
A non-trivial invariant set
$\Lambda$ in general supports an infinite number of invariant measures
some of which may be hyperbolic and some of which may not. At this point in
the theory it is not completely clear how one could distinguish situations
in which all invariant measure are hyperbolic and situations in which they
are not. For the moment we suggest the following definition: 
we say that a compact invariant set $\Lambda$ is 
\emph{totally hyperbolic} if all invariant 
measures with support on $\Lambda$ are hyperbolic. We conjecture
that the presence of (singular) invariant measures with 
zero Lyapunov exponent (a natural generalization of the 
indifferent fixed point case),
could be the main mechanism for slowing down of the mixing 
process and thus
giving rise to only subexponential rates of decay of correlations.

\begin{conjecture}\label{hyp}
The map \( f : \Lambda \to \Lambda \) exhibits
exponential decay of correlations
(with respect to every physical measure \(\mu\) with support in 
\( \Lambda  \)) if and only if $\Lambda$ is totally hyperbolic.
 \end{conjecture}
Conjecture~\ref{hyp} is true in the case of unimodal interval maps $f$ with
negative 
Schwarzian derivative. Indeed, as was shown in \cite{NS98},
$f$ has exponential decay of correlations if and only if
\[
\lambda_{per} := \inf \{ \frac1n \log |(f^n)'(p)| ;
n \geq 1, p \mbox{ is $n$-periodic}\,\} > 0,
\]
and \cite[Proposition 3.1]{BK} states that the Lyapunov
exponent of any $f$-invariant measure is at least $\lambda_{per}$.

Different degrees of hyperbolicity might also influence the effect of small
perturbations.
Tsujii \cite{Tsu93} showed that for generic one-parameter families
 unimodal maps satisfying a strong form of the Benedicks-Carleson conditions
 (and thus with exponential decay of correlations) are Lebesgue
 density points of similar maps.

 \begin{conjecture}
 For generic one-parameter families, maps with exponential decay of
correlations are
 Lebesgue density points of other maps with uniform exponential rates
 of decay of correlations.
  Maps with at least polynomial decay are Lebesgue density points of
 maps with (arbitrarily small) exponential decay.
 \end{conjecture}

 \subsection{Overview of the paper}

 Our strategy is to define a Markov return map
 \( \hat f= f^{R}:\Omega_{0}\to\Omega_{0} \) on a suitable
 neighbourhood of one of the critical points. We shall obtain
 estimates on the tail \( |\{x\in\Omega_{0}: R > n\}| \)
 of the return times and apply
 the general framework of L.-S. Young \cite{Y} linking these
 estimates with bounds for the decay of
 correlation. The general philosophy is that intervals
 outside a neighbourhood \( \Delta \) of the critical set grow
 exponentially fast (by a classical result of Misiurewicz \cite{mis},
see \cite{Man} for a $C^2$ version) and therefore for these intervals
during these times, the tail of the
return times decays exponentially fast.  Many intervals
however fall into \( \Delta \) before a good return to \( \Omega_{0}
\) occurs and are strongly contracted in the next iterate. We estimate
the time it takes for them to recover their original size
in terms of the derivatives along the critical orbits, which
in turn provides bounds for the decay of the
tail of the return time function.

In Section~\ref{Inducing to small scales} we consider intervals in \(
\Delta \) and use a binding argument to obtain estimates for their
growth in terms of the derivative along the appropriate critical orbit.
As mentioned above, similar arguments have been applied before, notably by
Jakobson \cite{Jak} and
Benedicks \& Carleson \cite{BC},
under stronger conditions on \( D_{n} \) and  on
the recurrence of the critical orbit. Here we have generalized the
argument to deal with  slow derivative
growth rates along the critical orbits and
arbitrary recurrence patterns.

In Section~\ref{large scales} we consider an arbitrary interval \(
J\subset I \) and show that there exists a partition \( \hat{\mathcal P}
\) of
\( J \) and a stopping time function \( \hat p \)
such that the images \( f^{\hat p(\omega)}(\omega) \) are uniformly
large for all \( \omega\in\hat{\mathcal P} \), i.e. almost every point
of \( J \) belongs to an interval which achieves \emph{large scale}.
We describe a combinatorial structure of \( f^{\hat p} \) on \( J \)
which keeps track of the pattern of returns to \( \Delta \) of each \(
\omega \). By combining this information with some
analytic estimates on the size of elements with given combinatorics, we
obtain key estimates on the size of the tail \( \{x\in J: \hat p > n\} \)
of the stopping time function \( \hat p \). A variety
of arguments is used here to deal with the various possible
rates (polynomial, stretched exponential or exponential).

In Section~\ref{The full return map} we show that once an interval
has achieved large scale there is a fixed proportion of it which
has a full return to the
original interval \( \Omega_0 \) within a fixed number of iterates.
It follows that the transition from
large scale to full return occurs exponentially fast and does not
significantly affect the tail estimates.
We also state precisely the results of Young which we apply
to our return map to obtain the conclusions of our theorems.\\
\\
{\bf Acknowledgements:}
The research for this paper was partly supported by the
PRODYN program of the European Science Foundation. S. Luzzatto also
acknowledges the financial support of EPSRC grant No. GR/K86329.

The suggestions of the referee have considerably improved
the present exposition. We acknowledge them gratefully.

\section{Inducing to small scales}
\label{Inducing to small scales}

We define a partition \( \mathcal P \) of a critical neighbourhood
\( \Delta \) and a stopping time function \( p \) such that the
induced map \( F=f^{p} \) on \( \Delta \) is expanding. The images of
partition elements are not uniformly large, i.e. \(
\inf\{|f^{p(\omega)}(\omega)| : \omega\in\mathcal P\}=0 \), and
therefore we call this \emph{inducing to small scales}.

\subsection{Definitions and notation}
\label{small scales map}

\begin{lemma}\label{equiv}
The conditions ($*$) and ($**$) are equivalent.
\end{lemma}

\begin{proof}
Condition ($*$) implies ($**$) because if we take $\gamma_n$ so that
$\gamma^{2\ell-1}=D_n^{-1}$ then
$[\gamma_n^{\ell-1} D_n(c)]^{-1/\ell}=\gamma_n=
D_n ^{-1/(2\ell-1)}$,
so the terms in each of the two sums in ($**$) are equal to each
other and equal to those in ($*$).

To see that ($**$) implies ($*$) note that by the duality
of $l^p$ and $l^q$ when $\frac{1}{p}+\frac{1}{q}=1$,\,\,
$\sum a_n^p<\infty$, $\sum b_n^q<\infty$ implies
$\sum a_n b_n < \infty$. Assume ($**$) holds and
take $a_n^p=\gamma_n$, $b_n^q=[\gamma_n^{\ell-1} D_n(c)]^{-1/\ell}$,
$1/q=\ell/(2\ell-1)$ and therefore $1/p=1-1/q=(\ell-1)/(2\ell-1)$.
Then $\sum a_n^p$ and $\sum b_n^q$ are both finite
and therefore $\sum a_n b_n < \infty$. But since
$a_n b_n=\gamma_n^{\frac{1}{p}}\gamma_n^{-\frac{\ell-1}{\ell}\frac{1}{q}}
D_n^{\frac{-1}{\ell q}} = D^{\frac{1}{2\ell-1}}$,
condition ($*$) follows.
\end{proof}

We use the symbol $\approx$ to indicate that two terms are equal
up to a factor depending only on $f$.
Because $\ell < \infty$,
$$
|f'(x)| \approx
|x-c|^{\ell-1}.
$$
for all $x \in X$ close to $c$.
Also there exists $\tau$ such that
\begin{equation}\label{estimate Dfx}
\frac{|f'(x) - f'(y)|}{|f'(x)|} \leq \tau
\frac{|x-y|}{|x - \Crit|},
\end{equation}
for all $x,y$ such that
$|x-y| \leq \frac12 \max\{ |x-\Crit|, |y-\Crit|\}$.
Here $|x- \Crit| = \min \{ |x-c|; c \in \Crit \}$.
Let
$\Gamma = \exp(\tau \sum_{j=1}^{\infty} \frac{\gamma_j}{1-\gamma_j})$.

For $x \in X$, let $c = c(x) \in \Crit$ be the critical point closest to
$x$. This is well defined for $x$ sufficiently close to $\Crit$.
Given a critical neighbourhood $\Delta$ of $\Crit$ we define the
\emph{binding period} as follows:
If $x \in \Delta$, then
\begin{equation}\label{bind}
 p(x) :=
  \max\{p: |f^k(x)-f^k(c)| \leq \gamma_{k} |f^k(c)-\Crit|
 \quad\forall\  k\leq p - 1 \},
  \end{equation}
while $p(x) := 0$ if $x \notin \Delta$.
Clearly $p \to \infty$ monotonically as $x \to c(x)$.
In order to choose the size of our critical neighbourhood \( \Delta \)
we need the following lemma.
\begin{lemma}\label{49}
Suppose that $G_p \geq 0$ and $\sum_p G_p < \infty$.
Then for any $\zeta > 0$ there exists $p_0$ such that
$$
P = \sum_{s \geq 1} \sum_{
\substack{(p_1, \dots, p_s)\\
p_i \geq p_0 }}
\prod_{p_i} \zeta G_{p_i} \leq 1.
$$
\end{lemma}

\begin{proof}
Let $S_0 = \sum_{p \geq p_0} \zeta G_p$. Then both $S_0$ and
$S := \sum_{s \geq 1} S_0^s$ tend to $0$ as $p_0 \to \infty$.
Developing term by term we see that
$P \leq S$. This proves the lemma.
\end{proof}

\begin{lemma}\label{LemBBC}
There exists $\kappa > 0$ such that
for all $\delta_0 > 0$, there exists $\delta \in (0,\delta_0)$
such that for $\Delta = \cup_c (c-\delta, c+\delta)$
and every $x$
\begin{equation*}\tag{BBC}
|(f^n)'(x)| \geq \kappa \text{ for } n = \min\{ i \geq 0; f^i(x) \in
\Delta\}.
\end{equation*}
\end{lemma}
We call this property {\em bounded backward contraction}. 
In an earlier version of this paper \cite{BLS}, we had to state (BBC)
as an assumption. For (symmetric)
S-unimodal maps, (BBC) is well-known to hold, cf. \cite{Guc},
and recently the multimodal case
it is proven in \cite{BvS1}. 
It is essential for (BBC) that all critical orders are
the same, see the counterexamples in \cite[Section 5]{BvS1}.

Taking advantage of Lemma~\ref{49} and condition ($**$) we
fix for the rest of the paper a critical neighbourhood
\( \Delta = \Delta_{\delta} = \cup_{c \in \Crit} (c-\delta,c+\delta) \)
where \( \delta >0 \) is such that
(BBC) holds and so small that
\begin{equation}\label{p delta}
\sum_{s \leq n} \sum_{
\substack{(p_1, \dots, p_s)\\
\sum_i p_i \leq n \\
p_i \geq p_{\delta} }}
\prod_{p_i} \zeta (\gamma_{p_i}^{\ell-1} D_{p_i}(c))^{1/\ell}
\leq 1
\end{equation}
for all $c \in \Crit$,
$p_{\delta} := p(c \pm \delta)$, $\zeta = 64 K_0/\kappa C_0$,
$C_0$ the constant introduced in Lemma~\ref{bound expansion},
and $K_0$ a fixed Koebe distortion constant, which turns
out to be $\leq 16$.

For  \( p\geq 0 \) we let $I_p= \{x: p(x) = p\}$
denote the level sets of the function \( p \). 
Let \( \mathcal P \)
denote the corresponding partition of $X$. 
Note that since $p(x) \equiv 0$ outside $\Delta$, 
$I_0 = X \setminus \Delta$ is the ``zeroth'' partition element.
Notice  that \( I_{p} \) can be empty for some values of \( p \), and that
it has at most $2\# \Crit$ components.
Define \( F: X \to X \) by letting \( F(x) = f^{p(x)}(x) \) for \(
x \in \Delta \) and \( F(x)=f(x) \) for \( x\in X\setminus\Delta \).

\subsection{Expansion estimates}
\label{small scales expansion}

We have two main expansion estimates.

 \begin{lemma}[Derivative growth for pieces of orbit outside
 \protect\( \Delta \protect\)]
 \label{outside delta}
\lineclear
 There exist constants $C_{\delta} > 0$
and $\lambda_{\delta} > 0$ such that for every piece of
 orbit \( \{f^{i}(x)\}_{i=0}^{k-1} \) lying completely outside
 \( \Delta \)
we have
 \[ |(f^{k})'(x)| \geq C_{\delta} e^{\lambda_{\delta} k}.
  \]
If moreover \( f^{k}(x) \in \Delta\), then
  \[
  |(f^{k})'(x)| \geq \max \{\kappa, C_{\delta} e^{\lambda_{\delta} k}\}.
   \]
 \end{lemma}
 Notice that the first estimate clearly implies the second if
 \( k \) is large. The second however is extremely useful
 when considering small values of $k$.

 \begin{proof}
 The first estimate is well known for maps with negative Schwarzian
 derivative, and also for maps without periodic attractors
or neutral orbits (Ma\~n\'e's result), for see Chapter II in \cite{MS}.
So this covers our case.  The second statement follows from (BBC). 
 \end{proof}
The following expansion bound will be of importance. Let
$$
F'_{p}(c) := \min\{ |(f^p)'(x)| ;
x \in I_{p} \cap (c-\delta, c+\delta) \}.
$$
\begin{lemma}[Derivative growth for pieces of orbit starting in
\protect\( \Delta \protect\)]
\label{bound expansion}
\lineclear
There exists $C_0 > 0$ (independently of $\delta$ and hence $\Delta$)
such that for every $c \in \Crit$ and
\( p\geq p_{\delta} \) with \( I_{p}\neq \emptyset \) we have
\begin{equation}\label{F'}
F'_{p}(c) \geq C_0 [\gamma_p^{\ell-1} D_p(c)]^{1/\ell}.
\end{equation}
\end{lemma}
In the sequel, we will write $F'_p$ instead of $F'_p(c)$ when no
confusion can arise.
We shall need an intermediate result for the proof.

\begin{lemma}\label{distortion in binding}
For \( x\in\Delta \)  we have
  \[
\frac{|(f^{i})'(y)|}{|(f^{i})'(z)|}\leq \Gamma
\text{ for all } \ y, z\in [f(x), f(c)] \text{ and  all }
\ i \leq p (x)-1.
   \]
\end{lemma}

{\em Remark:} In Subsection~\ref{combinatorial structure} 
we will use this estimate on a slightly bigger 
interval than $[f(x), f(c)]$, but this does not seriously affect the
estimates. 

\begin{proof}
    Letting \( y_{j}=f^{j}(y) \) and \( z_{j}=f^{j}(z) \) for \(
    j\geq 0 \) we have
by the chain rule
 \[
 \left|\frac{(f^{i})'(y)}{(f^{i})'(z)}\right|=
 \prod_{j=0}^{i-1}
 \left|\frac{f'(y_{j})}{f'(z_{j})}\right|
 =
\prod_{j=0}^{i-1}
 \left(1+\frac{|f'(z_{j})-f'(y_{j})|}{|f'(y_{j})|}
 \right).
  \]
By \eqref{estimate Dfx},
\(  |f'(z_j)-f'(y_j)|/|f'(y_j)| \le
\tau |z_j-y_j|/|y_j-\Crit| \)
and so, using the elementary fact that
\( \log (1+x) \leq x \) for all \( x>0 \) we get
\[
\log  \left|\frac{(f^{i})'(y)}{(f^{i})'(z)}\right|
\leq \sum_{j=0}^{i-1} \log \left(1+
\tau  \frac{|z_{j}-y_{j}|}{|y_{j}-\Crit|}\right)
\leq \tau  \sum_{j=0}^{i-1}
\left(\frac{|z_{j}-y_{j}|}{|y_{j}-\Crit|}\right).
 \]
By definition of \( p \)
we have
\( |z_{j}-y_{j}|\leq |f^{j+1}(x)-f^{j+1}(c)|\leq
\gamma_{j+1}|f^{j+1}(c)-\Crit| \) and \( |y_{j}-\Crit| \geq
(1-\gamma_{j+1})|f^{j+1}(c)-\Crit| \).
Here $c$ is again the critical point closest to $x$.
Substituting these inequalities into
the last formula yields the desired statement.
\end{proof}

\begin{proof}[Proof of Lemma~\ref{bound expansion}]
Let $c$ be the critical point closest to $x$.
By Lemma~\ref{distortion in binding} and the fact that
\(
|f'(x)|\approx  |x-c|^{\ell-1} \approx
|f(x) -f(c)|^{(\ell-1)/\ell}
\)
we have
\begin{equation}\label{fp1}
|(f^{p})'(x)|\geq \frac{|f'(x)| D_{p-1}(c)}{\Gamma}
\approx
\frac{ |f(x)-f(c)|^{(\ell-1)/\ell} D_{p-1}(c)}{\Gamma}.
\end{equation}
By the Mean Value Theorem, the definition of \( p \) and the
distortion estimate in
Lemma~\ref{distortion in binding} we have
\[
\Gamma D_{p-1} |f(x)- f(c)|
\geq |f^p(x)- f^p(c)| \geq
\gamma_{p}| f^{p}(c)-\Crit|
\]
and therefore
\begin{equation}\label{fp2}
|f(x)- f(c)|
\geq \frac{\gamma_{p} |f^{p}(c)-\Crit|}{\Gamma D_{p-1}(c)}.
\end{equation}
Substituting  \eqref{fp2} into \eqref{fp1} gives
$$
|(f^{p})'(x)| \geq {\mathcal O}(\Gamma^{-2+1/\ell})
\gamma_p^{(\ell-1)/\ell} D_{p-1}(c)^{1/\ell}
|f^{p}(c)-\Crit|^{(\ell-1)/\ell}.
$$
Because $|f'(f^p(c))| \approx |f^p(c) - \Crit|^{\ell-1}$, the
Chain Rule gives
$$
|(f^{p})'(x)| \geq {\mathcal O}(\Gamma^{-2+1/\ell})
\gamma_p^{(\ell-1)/\ell} D_p(c)^{1/\ell}.
$$
This proves the lemma.
\end{proof}

\section{Inducing to large scales}
\label{large scales}

The main result of this section is the following

\begin{proposition}\label{tail of hat p}
Suppose that \( f \)
satisfies ($*$). Then there exist $\delta' > 0$
such that for all $\delta''>0 \) the following properties hold.  For an
arbitrary interval \( J \subset X \) with \( |J|\geq \delta'' \) there
exists a partition \( \hat{\mathcal P} \) of \( J \) (mod 0) and a
stopping time function \( \hat p: \hat{\mathcal P} \to \mathbb N \)
such that for all \( \omega\in\hat{\mathcal P} \),
\( \hat F |\omega := f^{\hat p(\omega)}|\omega \) is a
diffeomorphism with
uniformly bounded distortion
and $|\hat F (\omega)|= |f^{\hat p(\omega)}(\omega)| \geq \delta'$.
Moreover the following estimates hold:
\begin{description}
    \item[Summable case] Under no conditions on $d_n(c)$ 
other than which stem from ($*$)
    \[
    \sum_{n} |\{\hat p > n | J \}| < \infty.
    \]
\item [Polynomial case]
If $d_n(c) \leq C n^{-\alpha}$
for all $c \in \Crit$ and $n \geq 1$, then there exists $\hat C > 0$
such that
 $$
|\{ \hat p > n | J\}| \leq \hat C n^{-\alpha}.
 $$
\item [Stretched exponential case]
If $b_n(c) \leq C e^{-\beta n^{\alpha}}, \ \alpha\in
(0,1), \
\beta > 0$
for all $c \in \Crit$ and $n \geq 1$,
then for each $\hat \alpha\in (0,\alpha)$
there exist $\hat \beta, \hat C > 0$ such that
$$|\{ \hat p > n | J\}| \leq \hat C e^{-\hat \beta n^{\hat \alpha}}.$$
\item [Exponential case]
If $b_n(c) \leq C e^{-\beta n}$, $\beta >0$ for all $c \in \Crit$
and $n \geq 1$, then
there exist $\hat \beta, \hat C > 0$ such that
$$|\{ \hat p > n | J \}| \leq \hat C e^{-\hat \beta n}.$$
\end{description}
\end{proposition}

Let us try to clarify the role of the constants in this
proposition, and their interdependence.
In the previous section we have fixed $\delta$. By the
Contraction Principle (see e.g. \cite[Section IV.5]{MS}), there
exists $\delta'$ such that for each component
$W$ of $\Delta \setminus \Crit$ and each $n \geq 0$,
$|f^n(W)| \geq \delta'$. This is the $\delta'$ of
Proposition~\ref{tail of hat p}.

The expression \( |\{\hat p> n | J\}| \) denotes
the conditional probability
\( |\{x \in J ; \hat p(x) > n\}| / |J| \).
In Section~\ref{combinatorial structure}
we define and describe the combinatorics of the
partition \( \hat{\mathcal P} \) of $J$
and the stopping time \( \hat p \).
In Section~\ref{Metric estimates} we prove some key estimates on the size
of an interval \( \omega_{p_{1}, \dots, p_{s}}\in\hat{\mathcal P} \)
with a given combinatorics. In Section~\ref{stopping time
estimates} we combine these with some counting arguments to obtain
estimates on \( \{\hat p > n | J \} \).
Note that the supremum of \( \{\hat p > n | J \} \), when taken over all
intervals $J$, will be infinite, because tiny intervals take a long time 
to reach large scale. When applying 
Proposition~\ref{tail of hat p} in Section~\ref{The full return map}, 
we will fix the minimal interval length $\delta'' := 
\min\{ \delta'/3, |\Omega_0| \}$, where
$\Omega_0$ is an interval specified in Subsection~\ref{full returns}.
In this way, we obtain a bound
of $\{ |\{\hat p > n | J \}| ; |J| \geq \delta''\}$ which depends
only on $\delta''$, that is: the $\hat C$'s in Proposition~\ref{tail of hat p}
depend on $\delta''$ but not on $J$.


\subsection{Combinatorial structure of \protect\(\hat F\protect\)}
\label{combinatorial structure}

We start with any interval $\omega \in {\mathcal P}|J$, i.e.
$\omega = I_p \cap J$ for some $p \geq 0$,
and let
$\nu_1 = \min\{ n \geq 0; f^n(\omega) \cap \Delta \neq \emptyset\}$
be the first visit to $\Delta$.
Write $\tilde\omega = f^{\nu_1}(\omega)$.
There are two (mutually exclusive) cases:
\begin{itemize}
\item $|\tilde \omega| < \delta'$. We partition $\tilde \omega$ by
intersecting $\tilde \omega$ with the elements $\{ I_p\}$.
Each interval $I_p \cap \tilde \omega$ for $p > 0$ is labeled as
{\it Deep Return}. The interval $I_0 \cap \tilde \omega$ is taken together
with the interval $I_p \cap \tilde \omega$ adjacent to it.
Because $\tilde \omega$ is not too large compared to this particular
component of $I_p$, the estimates of the binding period of $I_p$
go through, see Lemma~\ref{distortion in binding} 
and the remark below it.

\item $|\tilde \omega| \geq \delta'$.
We cut off the outmost intervals of length $\epsilon$
from $\tilde \omega$, and stop with the remaining middle part;
it has reached {\it large scale}. The subinterval $\omega_0 \subset \omega$
such that $f^{\nu_1}(\omega_0)$ equals this middle part of $\tilde \omega$
is added to the partition $\hat {\mathcal P}$.
Here $\epsilon \ll \delta'$ is a constant to be fixed in
the proof of Lemma~\ref{P''}. The $\epsilon$ will be used effectively
in Lemma~\ref{fixedprop}.
For the moment it suffices to know that $\epsilon$ is smaller than
each component of $\Delta$ and each component of $X \setminus \Delta$.

The outmost intervals $\tilde \omega_{\pm}$ of length $\epsilon$
are partitioned by
intersecting them with the elements $\{ I_p\}$.
Each interval $I_p \cap \tilde \omega_{\pm}$ for $p > 0$ is labeled as
{\it Deep Return}. The interval $I_0 \cap \tilde \omega_{\pm}$
(if it exists) is labeled as {\it Shallow Return}.
\end{itemize}
Note that if $x \in \omega$, then $f^n(x) \in \Delta$ only if $n$
is in a binding period
of $x$, $x$ has a deep return or if $n = \hat p(x)$.
At shallow return times $n$, $f^n(x) \notin \Delta$.

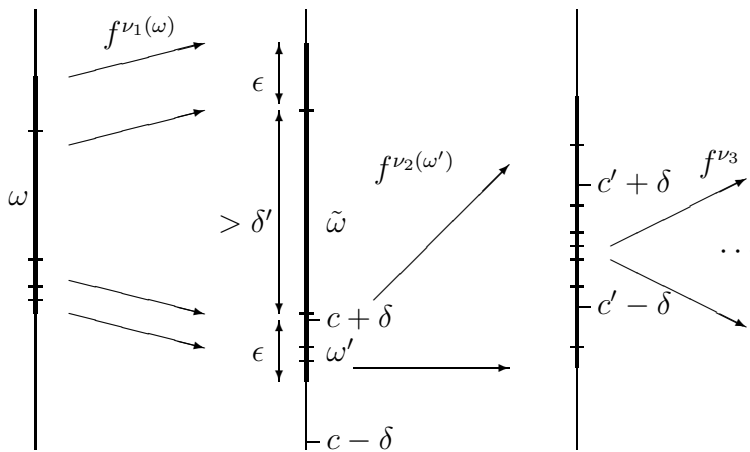
\begin{figure}
\unitlength=9mm
\begin{picture}(12,6.7)(0,0) \let\ts\mboxstyle

\thinlines
\put(1,0){\line(0,1){6.5}} \put(0.6,3.6){$\omega$}
\put(1.02,2){\line(0,1){3.5}} \put(0.98,2){\line(0,1){3.5}}
\put(5,0){\line(0,1){6.5}} 
\put(5.02,1){\line(0,1){5}} \put(4.98,1){\line(0,1){5}} 
\put(5.3, 3.2){$\tilde \omega$}  
\put(5.3, 1.3){$\omega'$}
\put(9,0){\line(0,1){6.5}}
\put(9.02,1.2){\line(0,1){4}} \put(8.98,1.2){\line(0,1){4}}
\put(0.9,2.8){\line(1,0){0.2}} \put(0.9,4.7){\line(1,0){0.2}}
\put(0.9,2.2){\line(1,0){0.2}} \put(0.9,2.4){\line(1,0){0.2}}

\put(4.9,2){\line(1,0){0.2}} \put(4.9,5){\line(1,0){0.2}} 
\put(4.6,2.2){\vector(0,1){2.8}} \put(4.6,4.8){\vector(0,-1){2.8}} 
\put(3.7,3.2){$>\delta'$}
\put(4.6,5.2){\vector(0,1){0.8}} \put(4.6,5.9){\vector(0,-1){0.8}} 
\put(4.2,5.3){$\epsilon$}
\put(4.6,1.1){\vector(0,1){0.8}} \put(4.6,1.8){\vector(0,-1){0.8}} 
\put(4.2,1.3){$\epsilon$}

\put(4.9,1.3){\line(1,0){0.2}} \put(4.9,1.5){\line(1,0){0.2}} 

\put(5.0,0.1){\line(1,0){0.2}} \put(5.3,0){$c-\delta$}
\put(5.0,1.9){\line(1,0){0.2}} \put(5.3,1.8){$c+\delta$}

\put(1.5,5.5){\vector(4,1){2}} \put(2,6){$f^{\nu_1(\omega)}$} 
\put(1.5,4.5){\vector(4,1){2}} 
\put(1.5,2){\vector(4,-1){2}} 
\put(1.5,2.5){\vector(4,-1){2}} 

\put(6,2.2){\vector(1,1){2}} 
\put(6,4){$f^{\nu_2(\omega')}$} 
\put(5.7,1.2){\vector(1,0){2.3}} 

\put(8.9,3){\line(1,0){0.2}} \put(8.9,3.2){\line(1,0){0.2}}
\put(8.9,2.8){\line(1,0){0.2}} \put(8.9,1.5){\line(1,0){0.2}}
\put(8.9,3.6){\line(1,0){0.2}} \put(8.9,4.5){\line(1,0){0.2}}
\put(8.9,2.4){\line(1,0){0.2}}

\put(9.0,2.1){\line(1,0){0.2}} \put(9.3,2.0){$c'-\delta$}
\put(9.0,3.9){\line(1,0){0.2}} \put(9.3,3.8){$c'+\delta$}

\put(9.5,3){\vector(2,1){2}} 
\put(10.8,4.1){$f^{\nu_3}$} 
\put(9.5,2.8){\vector(2,-1){2}} 
\put(11.1,2.8){$\cdots$} 

\end{picture}
\label{fig}
\caption{Construction of partitions $\hat {\mathcal P}_n$.
The middle part of $\omega$ reaches large scale after $\nu_1$ iterates;
the upper part has a shallow return, and the lower part has a deep return.
}
\end{figure}

Now let $\omega'$ be an interval which results from
this partitioning of $\tilde \omega$, which has not reached
large scale, see Figure~\ref{fig}. We first apply the binding period,
i.e. we take $f^p(\omega')$ for the stopping time $p = p(\omega')$
(which is possibly $0$, namely if $\omega' \cap \Delta = \emptyset$),
and then take the second return
$\nu_2 = \min\{ n \geq p(\omega');
f^n(\omega') \cap \Delta \neq \emptyset\}$.
Subdivide $f^{\nu_2}(\omega')$ according to the above rules, distinguishing
between large and deep returns.

Let $\hat{\mathcal P}_n$ be the partition which we obtain
by only considering at most $n$ iterates of $f$
and $\hat{\mathcal P}$ the partition of $J$
by considering all iterates of $f$.
We should emphasize that the procedure and hence the partitions 
depend on the choice made for $J$ and on $\epsilon)$.
For example, if $J_1$ and $J_2$ are two intersecting intervals,
then one could get two different partitions created at
a point $x\in J_1\cap J_2$.

Next we define the {\em stopping time at large scale} $\hat p_J$:
At points $x\in \omega$ where the procedure
eventually stops, i.e., for which there
exists $n > 0$ so that the $n$-th iterate
of the interval in $\hat{\mathcal P}_n$ containing $x$
has reached large scale, set $\hat p_J(x)=n$.
At other points $x\in J$ set $\hat p_J(x)=\infty$.
Finally to define $\hat F_J$, let
\( \hat J = \{x\in J: \hat p_J(x) <  \infty\}\)
and define \( \hat F_J: \hat J \to I \) by
$\hat F_J(x) = f^{\hat p_J(x)}(x)$.
We shall prove that
 \( \hat J = J  \) up to sets of zero Lebesgue measure.

Take $n<\infty$. To each $\omega_n \in \hat{\mathcal P}_n$ we
assign a formal {\em itinerary}
$$
(\nu_1, p_1), \dots, (\nu_s, p_s),
$$
consisting of the stopping times and lengths of the corresponding
binding periods;
$s$ is maximal for \( \nu_{s}\leq n \). Depending on the
depth of the return at time \( \nu_{s} \), \( \nu_{s}+p_{s}\)
can be arbitrarily large.  If the return at time $\nu_j$
is shallow, then $p_j = 0$.
If $\omega$ is an interval on which $(\nu_1, p_1), \dots, (\nu_{j-1},
p_{j-1})$
is constant and for which $\nu_j = \nu_j(\omega)$ is the next return to
$\Delta$, then the set $\{ x \in \omega; p_j(x) = p\}$ has at most
$4$ components. This maximum is attained when $|f^{\nu_j(\omega)}| \geq
\delta$, the radius of $\Delta$,
and the outmost intervals of size $\epsilon$ both contain a critical point.
It can happen that $f^{\nu_j}(\omega)$ covers many more critical
points, but since $\omega$ has reached large scale, $p_j$ is only defined
on the outmost intervals.
We will take care of this multiplicity in the estimates in
Subsection~\ref{stopping time estimates}.
But apart from this multiplicity, a sequence $p_1,\dots, p_s$
uniquely determines a partition element
$\omega_{p_1,\dots,p_s} \in \hat{\mathcal P}_n$
(or perhaps the empty set).
Indeed, $p_i$ determines the position of the $i$-th return
of $\omega_{p_1, \dots , p_s}$, and from the previous
$p_1, \dots, p_{i-1}$ and the starting interval $J$ one can
compute the next return time. Hence the information
$\nu_1, \dots, \nu_s$ is strictly speaking superfluous.
Observe however that there are many itineraries
that do not correspond to partition elements.
Note that $f^n$  is a diffeomorphism
on each interval from the partition
$\hat{\mathcal P}_n$.

For a given sequence $(p_1, \dots ,p_s)$,
let
$$S_d=\{i\le s; \nu_i \mbox{ is a deep return }\}=
\{i\le s; p_i>0\}$$
and
$$S_s =\{i\le s; \nu_i \mbox{ is a shallow return }\}=
\{i\le s; p_i=0\}=S\setminus S_d.$$
Moreover, let
$$
S_{s,s}=\{i<s; p_i = 0\mbox{ and }p_{i+1} = 0\}.
$$
Because each index in $S_s \setminus S_{s,s}$
either equals $s$
or is followed by an index in $S_d$, we get
\begin{equation}\label{S_s}
\# S_s \leq \# S_{s,s} + \#S_d + 1.
\end{equation}

\subsection{Metric and combinatorial estimates}
\label{Metric estimates}

\begin{lemma}\label{size}
Let $C = C_{\delta}$ and $\lambda = \lambda_{\delta}$
be as in Lemma~\ref{outside delta}. There exists \( K_{0}>0
\) independent of \( \epsilon \) and \( \rho\in (0,1) \)
($\rho \to 0$ as $\epsilon \to 0$), with the following properties.
For a given sequence $(\nu_{1}, p_1), \dots (\nu_{s},p_s)$ with \(
\nu_{s}\leq n \) we have
 \[
\frac{ |\omega_{p_1, \dots, p_s}| }{ |f^m(\omega_{p_1, \dots , p_s})| }
 \leq \min
 \{
 C^{-\# S_d} e^{ -\lambda(m-\sum_{i=0}^s p_i)},
 \left(\frac{K_0}{\kappa}\right)^{\# S_d} \! \! \rho^{\# S_{s,s}}
 \}
 \prod_{i\in S_d} (F'_{p_i})^{-1}
\]
for $m = \max\{ n, \nu_s + p_s\}$.
Moreover there exists \( T>0 \) which can be chosen arbitrarily large
if \( \epsilon \) is small, such that
$\nu_{i+1}-\nu_i\ge T$ whenever $p_i=p_{i+1}=0$.
\end{lemma}

 \begin{proof}
By construction,
\( f^m|\omega _{p_1,\dots, p_s } \)
is a diffeomorphism onto its image.
Take $x \in \omega_{p_1, \dots, p_s}$ and divide its orbit
into pieces separated by returns
(both deep and shallow):
$$
[1, \nu_1-1], \, [\nu_1, \nu_2-1], \, [\nu_2, \nu_3-1], \, \dots, \,
[\nu_s, m].
$$
Let $\nu_i<\nu_{i'}$ be two consecutive
deep returns. That is, assume that $p_i>0$,
$p_{i+1}=\dots=p_{i'-1}=0$ and $p_{i'}>0$
(with possibly $i'=i+1$).
Because each such interval lasts at least the corresponding
binding period, and in the
remaining time the point $x$ does not visit $\Delta$,
Lemma~\ref{outside delta} and the definition of \( F'_{p} \)
give
$$
|(f^{ \nu_{i'}- \nu_i})'(f^{\nu_i}(x))|
\geq
C e^{ \lambda(\nu_{i'}-(\nu_i + p_i)) } F'_{p_i}.
$$
Hence the chain rule and
the Mean Value Theorem show that
$$
|\omega_{p_1, \dots , p_s}|
\leq C^{-\# S_d} e^{-\lambda(m-\sum p_i)}
\, |f^m(\omega_{p_1, \dots , p_s})| \,
\prod_i 1/F'_{p_i}.
$$
To prove the other inequality,
let $\nu_i$ and $\nu_{i'}$ be subsequent deep return times.
First let us treat the step from $\nu_{i'-1}+p_{i'-1}$ to $\nu_{i'}$,
so assume that
$\nu_{i'-1} + p_{i'-1} < \nu_{i'}$.
(If $\nu_{i'-1} + p_{i'-1} = \nu_{i'}$ we can skip this step.)
Therefore
$f^{\nu_{i'-1} + p_{i'-1}}(\omega_{p_1, \dots, p_s}) \cap \Delta =
\emptyset$,
while there is at least one point $y \in \omega_{p_1, \dots, p_s}$
such that $f^{\nu_{i'}}(y) \in \Delta$. Lemma~\ref{LemBBC} yields that
$$
|( f^{ \nu_{i'} - (\nu_{i'-1} + p_{i'-1}) } )' (f^{\nu_{i'-1} +
p_{i'-1}}(y))|
\geq \kappa.
$$
Because $\nu_{i'} < n$,
$|f^{\nu_{i'}}(\omega_{p_1, \dots, p_s})| < \delta'$.
Take $H \supset f^{\nu_{i'-1}+p_{i'-1}}(\omega_{p_1, \dots, p_s})$
the largest interval on which
$f^{\nu_{i'} - (\nu_{i'-1}+p_{i'-1})}$ is monotone.
Then by the choice of $\delta'$,
$f^{\nu_{i'} - (\nu_{i'-1}+p_{i'-1})}(H)$ contains a $\frac13$-scaled
neighbourhood of $f^{\nu_{i'}}(\omega_{p_1, \dots, p_s})$.
Therefore the derivative of $f^{\nu_{i'} - (\nu_{i'-1}+p_{i'-1})}$
has distortion bounded by some $K_0 = K_0(\frac13) \leq 16$.
This follows from the Koebe Lemma, see \cite[Chapter IV]{MS}.
Hence
$$
|( f^{ \nu_{i'} - (\nu_{i'-1} + p_{i'-1}) } )' (f^{\nu_{i'-1} +
p_{i'-1}}(x))|
\geq \kappa/K_0.
$$

If $i' \leq i+2$, the same argument gives
$|( f^{ \nu_{i'} - (\nu_i + p_i) } )' (f^{\nu_{i'-1} + p_{i'-1}}(x))|
\geq \kappa/K_0$, and indeed in this case there are no
entries of $S_{s,s}$ between $i$ and $i'$.

If $i' > i+2$, then the differences
$\nu_{i+2}-(\nu_i+p_i)$, $\nu_{i+3} - \nu_{i+2}$,\dots,
$\nu_{i'-1} - \nu_{i'-2}$ are all large if $\epsilon$ is small.
Indeed, in these times an interval of size $\epsilon$ must have expanded
to an interval of size $\delta' \gg \epsilon$.
Because $x$ does not visit $\Delta$ during these iterates
(recall that the binding periods at shallow returns have length
$0$), the first part of Lemma~\ref{outside delta} gives
$$
|( f^{\nu_{i'-1} - (\nu_i + p_i)} )' (f^{\nu_i+p_i}(x)) |
\geq \left(\frac{1}{\rho}\right)^{i'-(i+2)},
$$
where $\rho \to 0$ as $\epsilon \to 0$.
(In Lemma~\ref{P''} we will fix $\rho$ at $1/8$.)
Adding the numbers $i'-i-2$ (running over all pairs 
$(i',i)$ of subsequent deep returns
with $i' > i+2$) gives $\# S_{s,s}$. This proves the lemma.
\end{proof}

\begin{lemma}\label{distortion hat p}
There exists $K > 0$ depending only on $\epsilon$ such
that for all starting intervals $J$ and
$\omega \in \hat{\mathcal P}$,
the distortion of $\hat F_J|\omega$ is bounded by $K$.
\end{lemma}

\begin{proof}
If $\hat p(\omega) = n$, then by construction, there is an interval
$T \supset \omega$ such that $f^n$ maps $T$ monotonically
onto an $\epsilon$-scaled neighbourhood of $f^n(\omega)$,
i.e. both components of $f^n(T) \setminus f^n(\omega)$ have size
$\geq \epsilon |f^n(\omega)|$.
The Koebe Principle (see e.g. \cite{MS} and in the setting
when we do not assume $Sf<0$, \cite{vSV}) gives the result.
\end{proof}
The following lemma contains combinatorial estimates needed
in the next section.

\begin{lemma}\label{Stirling}
Let $N_{k,s}$ be the number of integer
sequences $(p_1, \dots, , p_s)$ such that
$p_1 + \dots + p_s = k$ and $p_i \geq 0$ for all $i$.
Let $N^+_{k,s}$ be the same number for sequences with
$p_i > 0$ for all $i$.
Then 
\begin{equation}\label{Nks}
N_{k,s} \leq 2^s \max_{j\le k} N^+_{k,j} < 2^{k+s-1}.
\end{equation}
Given $\zeta > 0$ small and $\alpha \in (0,1]$, there exists
$\hat \zeta = \hat \zeta(\zeta,\alpha)$ with $\hat \zeta \to 0$ as $\zeta \to 0$
such that
\begin{equation}\label{N+ks}
N^+_{k,s} \leq \left\{
\begin{array}{ll}
e^{ \hat \zeta k^{\alpha} \log k} &\text{ if } s \leq \zeta k^{\alpha}, \\
e^{\hat \zeta k} &\text{ if } s \leq \zeta k. \\
\end{array} \right.
\end{equation}
\end{lemma}

\begin{proof}
Say $j \geq 1$ terms in the sum $p_1 + \dots + p_s$ are nonzero.
There are $\binom{s}{j}$
ways to distribute them. This gives
$$
N_{k,s} = \sum_{j=1}^s \binom{s}{j}
N^+_{k,j}
\leq \max_j N^+_{k,j} \sum_{j=1}^s \binom{s}{j}
= (2^s-1) \max_j N^+_{k,j}.
$$
Next, supposing that $p_i > 0$ for all $i$, there are
$s-1$ partial sums $\sum_{i=1}^j p_i$ different from $k$.
Therefore $N^+_{k,j} = \binom{k-1}{j-1}
$.
In particular, $\max_j N^+_{k,j} \leq 2^{k-1}$.
This proves \eqref{Nks}.
Let us estimate $N^+_{k,s}$ more precisely if $s \leq \zeta k^{\alpha} \leq
k/2$.
By Stirling's formula
$$
\binom{k}{s} \leq \frac{k^k}{(s^s)(k-s)^{k-s}}
\leq
\left(\frac{k}{k-s}\right)^k \cdot \left(\frac{k-s}{s}\right)^s
\leq (1+\frac{2s}{k})^k \cdot \left(\frac{k}{s}\right)^s.
$$
Because $(1+\frac{2s}{k})^k \leq \exp( k \log(1+\frac{2s}{k})) \leq
\exp(2s)$,
it follows that
$$
N^+_{k,s} < e^{2s} \left( \frac{ k^{1-\alpha} }{ \zeta } \right)^{\zeta
k^{\alpha}}
\leq e^{ \zeta(2-\log \zeta + (1-\alpha) \log k)k^{\alpha} }
$$
if $s \leq \zeta k^{\alpha}$. If $s \leq \zeta k$, this simplifies to
$N^+_{k,s} < e^{\zeta (2-\log \zeta) k}$, proving \eqref{N+ks}.
\end{proof}

\subsection{Stopping time estimates}
\label{stopping time estimates}

The aim of this section is to estimate the
tail behaviour of the return time function \( \hat p \), i.e. to
obtain an upper bound for the Lebesgue measure of the
set \( \{x \in J: \hat p(x) >n \} \) which we shall
henceforth (suppressing the dependence on $J$) 
denote by \( \{\hat p > n\} \). 
We shall always assume the
notation of Proposition~\ref{tail of hat p}, particularly when
referring to the polynomial, stretched exponential and exponential
cases.

We fix \( n \) for the rest of this section.
For each \( \omega\in \hat{\mathcal P}  \)
we consider
the sequence \( p_1, \dots, p_s \) as defined above,
with some terms possibly equal to \( 0 \). Recall that \( s \) is
given by the number of returns occurring before time \( n \).
Let \( \eta > 0\) be
a small constant to be determined in Lemma~\ref{P'}.
The set of partition elements $\omega \in \hat{\mathcal P}_n$
with $\hat p|\omega > n$ can be divided into
\[
\hat{\mathcal P}'_n=\left\{\omega\in\hat{\mathcal P}_{n} ;\, \hat
p|\omega > n, \,
\sum_{i=1}^s p_i \leq \eta n \right\}
\]
and
\[
\hat{\mathcal P}''_n
=\left\{\omega\in\hat{\mathcal P}_{n}; \, \hat
p|\omega > n, \,
\sum_{i=1}^s p_i > \eta n \right\}.
 \]
 Clearly we have
 \[
 |\{\hat p > n \}|
 =\sum_{\omega\in\hat{\mathcal P}'_n}|\omega| +
 \sum_{\omega\in\hat{\mathcal P}''_n}|\omega|.
  \]
  To treat the exponential and stretched exponential case we shall
  need to subdivide \( \hat{\mathcal P}''_n \) further into
  \[
   \hat{\mathcal P}''_{n-}=\{\omega\in \hat{\mathcal P}''_{n} ;
   s \leq \rho n^{\hat \alpha}\} \quad\text{and}\quad
   \hat{\mathcal P}''_{n+}=\{\omega\in \hat{\mathcal P}''_{n} ;
   s > \rho n^{\hat \alpha}\},
  \]
where \( \hat \alpha \in (0,1] \) and $\rho > 0$ will be chosen  below.

 Intuitively elements in \( \hat{\mathcal P}'_n \)
 spend most of their time in the
 ``uniformly expanding'' region \( X\setminus\Delta \).
 Thus intervals are growing in size at a uniform exponential rate
 and achieve large scale exponentially fast.
 Elements of \( \hat{\mathcal P}''_n \) on
 the other hand  spend much time in binding periods.
 In this case the upper bound will more closely
reflect the expanding properties of the critical orbit.
We shall apply various combinations of the estimates obtained in
Section~\ref{Metric estimates} to obtain bounds on the total measure of
the elements of the subpartitions defined above under the required
assumptions on the growth of $D_n$. 

\begin{lemma}\label{P'}
For any \( \theta > 0 \) there exists \( \eta_0 > 0 \)
such that for all \( 0<\eta<\eta_0\) and
for all $n$ sufficiently large,
\[
\sum_{\omega\in\hat{\mathcal P}'_n}|\omega|\leq
e^{- (\lambda - \theta) n}.
\]
 \end{lemma}
 \begin{proof}
As in Lemma~\ref{Stirling}, let \( N_{k,s} \) denote the number of possible
sequences
\( (p_1 , \dots, p_s) \) with \( p_i \geq 0 \) and
\( p_1+\dots+p_s=k \).
Then by the definition of \( \hat{\mathcal P}' \) and
the first statement of Lemma \ref{size} we have
\begin{equation}\label{henk}
\sum_{\omega\in\hat{\mathcal P}'}|\omega| \leq
\sum_{s=1}^{n} \sum_{k=0}^{\eta n}
\sum_{
\substack{(p_1, \dots, p_s)\\
\sum p_i = k}}
|\omega_{p_1, \dots, p_s}|
\leq \sum_{s=1}^n \sum_{k=0}^{\eta n}
4^s N_{k, s} C^{-s} e^{-\lambda(1 - \eta) n}.
\end{equation}
Here the factor $4^s$ expresses the maximal number of
components of $I_p$ for each return, see the argument in
Subsection~\ref{combinatorial structure}.
We use the bound $N_{k,s} < 2^{k+s}$ from Lemma~\ref{Stirling}.
Recall that $k \leq \eta n$. 
Since $\nu_{i+1}-\nu_i\ge T$ when $p_i=p_{i+1}=0$ (see the previous lemma)
formula~\eqref{S_s} gives $s = \#S_d + \#S_s \le 2\#S_d + \#S_{s,s} +1
\le 2\eta n + n/T + 1$. So in (\ref{henk}), $s$ only ranges up to this bound.
Writing $\eta' = (3 \eta + 1/T + 1/n)
(3\log 2 + \log C^{-1})$,
we get
$4^s N_{k,s} C^{-s} \leq e^{\eta' n}$.
Taking $\theta = 2 (\eta + \eta'/ \lambda)$
and substituting in \eqref{henk} gives
\[
\sum_{\omega\in\hat{\mathcal P}'_n}|\omega| \leq
n \sum_{k=1}^{\eta n} e^{\eta' n} e^{-\lambda(1-\eta) n}
\leq \eta n^{2} e^{-(\lambda - \frac{\theta}{2}) n} \leq
e^{-(\lambda-\theta) n}
 \]
provided $\eta$ and $\eta'$ are sufficiently small
and $n$ sufficiently large.
\end{proof}

\begin{lemma}\label{P''}
Recall from \eqref{dn} that $d_n(c) = 
\min_{i < n} (\gamma_i/D_i(c))^{1/\ell}|f^i(c) -
\Crit|$.
Fix $L \in \{ 1, \dots, n\}$ arbitrary and let
$$
\hat d_{n,s}(c) = d_i(c) \text{ for } 
i = \max\{ \lceil \frac{\eta n}{2s^2} \rceil , L\}.
$$
Write $s(\omega) = s$ if the itinerary $(p_1,\dots, p_s)$ of
$\omega$ has length $s$. 
For any \( \eta>0 \) there exists $C_1  > 0$
such that
$$
\sum_{\substack{ \omega\in\hat{\mathcal P}''_n \\ s(\omega) \geq L} }
|\omega| \leq
C_1 \max_{c \in \Crit} \sum_{s = L}^n 2^{-s}  \hat d_{n,s}(c).
$$
\end{lemma}

\begin{proof}
Given a sequence $(p_1,\dots,p_s)$, let $p_{j'}$ be the first term
such that $p_{j'} \geq \eta n/(2 {j'}^{2})$. Because $p_1 + \dots + p_s \geq
\eta n$, such $j'$ exists. Take $j = \max\{ L, j'\}$.

Let $\tilde \omega_{p_1,\dots, p_j}$ be the union of adjacent intervals
$\omega_{p_1, \dots, p_{j-1},p}$
with common return times $\nu_1, \dots, \nu_j$ and $p \geq p_j$.
Then $f^{\nu_j}$ maps $\tilde \omega_{p_1, \dots, p_j}$
diffeomorphically into an interval $(x,y)$ such that
$p(x) , p(y) \geq p_j$.
Assume without loss of generality that $|x-c| \geq |y-c|$.
Therefore, for each $i < p_j$,
\begin{align*}
\gamma_i |f^i(c) - \Crit|  & \geq |f^i(x) - f^i(c)| \\
&\geq \Gamma^{-1} D_{i-1}(c) |f(x) - f(c)| \\
&\geq {\mathcal O}(1/\Gamma) D_{i-1}(c) |x-c|^{\ell} \\
&\geq {\mathcal O}(1/\Gamma)
\frac{ D_i(c) |x-c|^{\ell} }{ |f^i(c) - \Crit|^{\ell-1} }.
\end{align*}
This gives
\begin{align*}
|x-y| \leq 2 |x-c| &\leq
{\mathcal O}(2\Gamma^{1/\ell}) \, d_p(c)  \\
&\leq
{\mathcal O}(2\Gamma^{1/\ell})
\max_{p \geq \eta n/2j^2} d_p(c)
= {\mathcal O}(2\Gamma^{1/\ell}) \, \hat d_{n,j}(c).
\end{align*}
Let $S_s$ and $S_d$ be the indices $\leq j$ corresponding to
shallow respectively deep returns. Also let $S'_d = S_d \setminus \{ j \}$
and let $S_{s,s}$ be the indices $\leq j$
of shallow returns that are followed by another shallow return.
Now Lemma~\ref{size} applied to $\tilde \omega_{p_1, \dots, p_j}$
and the iterate $\nu_j$ gives
\begin{align*}
\sum_{\substack{ \omega\in\hat{\mathcal P}''_n \\ s(\omega) \geq L} } |\omega|
&\leq \sum_{j = L}^n
 \sum_{(p_1, \dots, p_j)} |\tilde \omega_{p_1, \dots, p_j}| \\
&\leq
\sum_{j = L}^n
{\mathcal O}(2\Gamma^{1/\ell}) \max_{c \in \Crit} \hat d_{n,j}(c)
\sum_{(p_1, \dots, p_{j-1})} \!\!\!\!
4^j \left(\frac{K_0}{\kappa}\right)^{\# S_d} \rho^{\# S_{s,s}}
 \prod_{i \in S'_d} \frac{1}{F'_{p_i}}.
\end{align*}
The factor $4^j$ expresses the
different components of the level sets $I_p$ that intersect
forward iterates of $\omega$
(see the argument in Subsection~\ref{combinatorial structure}), 
and the factor
$(K_0/\kappa)^{\# S_d} \rho^{\# S_{s,s}}$ comes from Lemma~\ref{size}.
Using \eqref{S_s} and the fact that \( \#S'_{d}=\#S_{d}-1 \) we can write
\(  4^j = 2^{-j} 8^j \) and
\[
8^{j}= 8^{\#S_{s}+\#S_{d}}\leq 8^{\#S_{s,
s}+2\#S_{d}+1}=  8^{\#S_{s,s}} 64^{\#S_{d}} 8 =
512 \  8^{\#S_{s,s}} 64^{\#S'_{d}}.
\]
Then
$$
\sum_{(p_1, \dots, p_{j-1})}
\!\!\!\!
4^j \left(\frac{K_0}{\kappa}\right)^{\# S_d} \rho^{\# S_{s,s}}
\prod_{i \in S'_d} \frac{1}{F'_{p_i}}
\leq
2^{-j}
\frac{512 K_0}{\kappa}
\sum_{ (p_1, \dots, p_{j-1})}
\!\!\!\!
(8\rho)^{\# S_{s,s}} \prod_{i\in S'_d} \frac{64 K_0}{ \kappa F'_{p_i}}.
$$
Take $\epsilon$ in Lemma~\ref{size} so small that
$\rho = \frac18$ and
recall that $p_i \geq p_{\delta}$ for all $i \in S_d$.
Therefore Lemma~\ref{bound expansion} and formula~\eqref{p delta}
(with $\zeta = 64 K_0/\kappa$ ) give that
$$
\sum_{ (p_1, \dots, p_{j-1}) }
(8\rho)^{ \# S_{s,s} }
\prod_{i\in S'_d} \frac{ 64 K_0 }{ \kappa F'_{p_i} } \leq 1.
$$
By \eqref{F'}, the lemma follows with
$C_1 = {\mathcal O}(2\Gamma^{1/\ell})
512 K_0/\kappa$.
\end{proof}

The previous lemma is not so useful in the exponential and stretched
exponential cases for relatively small values of \( s \). Indeed, consider
for example the situation that $d_p = e^{-\beta p}$.
Then the term in the sum in Lemma~\ref{P''}
corresponding to $s=\sqrt n$ gives 
$C_1 2^{-\sqrt n} \cdot e^{-\eta \beta/2}$.
Clearly this decreases merely subexponentially in $n$. Let us improve on
this.

\begin{lemma}\label{P'''}
Assume that there exists $C,\beta>0$ and $\alpha\in (0,1]$
such that $D_n^{-1/\ell}\leq C e^{-\beta n^{\alpha}}$ for all
$n$.
Then for each $\hat \alpha \in (0,\alpha)$ (or $\hat \alpha = 1$ if $\alpha
= 1$)
there exists $\rho, C', \beta'>0$ such that
\[
\sum_{\omega\in\hat{\mathcal P}''_{n-}}|\omega| \leq
C' e^{-\beta' n^{\alpha}} \]
for all $n$.
Note that the set ${\mathcal P}''_{n-}$ depends on $\rho$ and $\hat \alpha$.
\end{lemma}

\begin{proof} First notice that since $\alpha\in (0,1]$
one has $p_1^\alpha+p_2^\alpha\ge (p_1+p_2)^\alpha$.
Using Lemmas~\ref{bound expansion} and \ref{size} this gives that 
there exists $\beta''>0$ and $C$ such that
$$
|\omega_{p_1, \dots, p_s}| \, \leq \,
C^{-s} \prod_{i=1}^s \frac1{F'_{p_i}} 
\, \leq \, 
C^{-s} C_0^{-s} \prod_{i=1}^s \max_{c \in \Crit} b_{p_i}(c) 
\, \leq \,
C^{-s} e^{-\beta'' (\sum p_i)^\alpha }.
$$
Reasoning as in the proof of Lemma~\ref{P'}, we write \( k=
p_{1}+\dots+p_{s} \) and we obtain
\begin{equation}\label{fipa}
\sum_{\omega \in \hat{\mathcal P}''_{n-} } |\omega|
\leq \sum_{s=1}^{\rho n^{\hat \alpha}} \sum_{k=\eta n}^\infty
4^s N_{k, s} C^{-s} e^{-\beta'' {k^\alpha} }.
\end{equation}
Taking $\zeta=\rho/\eta^{\hat \alpha}$ respectively
$\zeta=\rho/\eta$ in Lemma~\ref{Stirling}, we get that for some
$\hat \rho = \hat\rho(\rho, \eta,\hat\alpha)$ with
$\hat \rho \to 0$ as $\rho \to 0$, 
$$
N_{k,s} \leq \left\{
\begin{array}{ll}
2^s e^{ \hat \rho k^{\hat \alpha} \log k} &\text{ if } s \leq 
(\rho/\eta^{\hat \alpha}) k^{\hat \alpha}, \\
2^s e^{\hat \rho k} &\text{ if } s \leq (\rho/\eta) k. \\
\end{array} \right.
$$
(The second case applies when $\alpha=1$.)
Because $\hat \alpha \leq \alpha$ and taking $\rho$ and therefore $\hat \rho$
sufficiently small, we get
in either case
$$
\sum_{k = \eta n}^\infty 4^s N_{k,s} C^{-s} e^{-\beta'' {k^\alpha}}
\leq 8^s C^{-s} e^{-\beta'' (\eta n)^{\alpha}/2 }.
$$
Using again that $\hat \alpha \leq \alpha$ and the fact that
$\rho$ is small, inequality~\eqref{fipa} gives
\[
\sum_{\omega \in \hat{\mathcal P}''_{n-}} |\omega|
\leq
\sum_{s=1}^{\rho n^{\hat \alpha}} \sum_{k=\eta n}^{\infty}
4^s N_{k, s} C^{-s} e^{-\beta'' {k^\alpha}} 
\leq
C' e^{-\beta'' (\eta n)^{\alpha}/4 },
\]
for some constant $C'$.
This proves the lemma with $\beta' = \eta^{\alpha} \beta''/4$.
\end{proof}

\begin{proof}[Proof of Proposition~\ref{tail of hat p}]
 We show first of all that
 \( \hat J \) has full measure in \( J \), i.e. \( |\{\hat p > n\}|
 \to 0 \) as \( n\to \infty \).
 By (for example) Lemma~\ref{outside delta},
it follows that almost all $x\in J$,
$f^n(x)$ accumulates onto $\Crit$.
Hence $x$ has infinitely many deep return times, and
it is contained in sets of the form
$\omega_{p_1, \dots , p_s}$ for itineraries of arbitrary length $s$.
Because $\hat p(\omega_{p_1, \dots , p_s}) \geq s \to \infty$
as $s \to \infty$,
the proofs of Lemmas~\ref{P'} and \ref{P''} show that
$\sum_{(p_1, \dots, p_s)} |\omega_{p_1, \dots ,p_s}| \to 0$
as $s \to \infty$.
Therefore $| J \setminus \hat J| = 0$.

To prove the remaining estimates in the four cases mentioned in the
proposition, notice that we have exponential bounds for
\( \hat{\mathcal P}'_{n} \) and therefore we only need to concentrate
here on \( \hat{\mathcal P}''_{n} \).
The sequence $\{ \hat d_{n,s}(c) \}$ is decreasing in $n$,
and for each $k$ there are
at most $\#\{ n ; k-1 \leq \eta n/(2 s^{2}) \leq k\} \leq 2 s^{2}/\eta$
numbers $n$ such that $k = [\eta n/(2 s^{2})]$.
Therefore, using Lemma~\ref{P''} with $L = 1$:
\begin{align*}
\sum_{n \geq 1} \sum_{s=1}^n 2^{-s} \hat d_{n,s}(c)
&\leq
\sum_{s \geq 1}  \frac{2 s^{2}}{\eta} 2^{-s} \sum_{k \geq 1}
[\gamma_k/D_k(c)]^{1/\ell} \ |f^k(c)- \Crit| \\
&\leq \frac{12}{\eta} \sum_{k \geq 1} [\gamma_k / D_k(c)] ^{1/\ell} \\
&\leq  \frac{12}{\eta} \sum_{k \geq 1} [\gamma_k^{\ell-1} D_k(c)]
^{-1/\ell}.
\end{align*}
Hence the summable case follows from ($**$).
Lemma~\ref{P''} with $L = 1$ gives for the polynomial case
$$
\sum_{\omega \in {\mathcal P}''_n } |\omega|
\leq C_1 \max_{c \in \Crit} \sum_{s = 1}^n 2^{-s} \hat d_{n,s}(c)
\leq C_1 \sum_{s=1}^n 2^{-s} \left( \frac{2s^2}{\eta n} \right)^{\alpha}
\leq 12 C_1 \eta^{-\alpha} n^{-\alpha}
$$
as required.
In the exponential and stretched exponential cases we use
Lemma~\ref{P''} applied to
$\hat{\mathcal P}''_{n+}$ with $L = \rho n^{\hat \alpha}$ to get
$$
\sum_{\omega \in {\mathcal P}''_{n+} } |\omega|
\leq C_1 \max_{c \in \Crit}
\sum_{s \geq \rho n^{\hat \alpha}} 2^{-s} \hat d_{n,s}(c)
\leq C_2 e^{-(\log 2)\rho n^{\hat \alpha} }
$$
for some $C_2 > 0$. Lemma~\ref{P'''} takes care of the remaining
collection ${\mathcal P}''_{n-}$.
\end{proof}

\section{The full return map}\label{The full return map}

In this section we construct the full return map \( \hat f:\Omega_0 \to
\Omega_0 \)
and carry out its tail estimates.

\begin{proposition}
\label{return times}
    Suppose that \( f \) satisfies (\( * \)). Then
    for any \( c\in\Crit\cap X \) there exists a neighbourhood \( \Omega_{0}
    \) of \( c \), a countable partition \( \mathcal Q \) of \(
    \Omega_{0} \) (mod 0) and a return time function \( R: \mathcal Q
    \to \mathbb N \) with the following properties. For each \(
    \omega\in\mathcal Q \), \( \hat f := f^{R} \)
    maps \( \omega \) to \(
    \Omega_{0} \) diffeomorphically with bounded distortion:
    letting
\begin{equation}\label{separation}
s(x, y) = \min\{n ; \hat f^{n}(x), \hat f^{n}(y)
\text{ belong to different
elements of } \mathcal Q \},
\end{equation}
there exists $\beta \in (0,1)$ and $C > 0$
such that for all $\omega \in {\mathcal Q}$ and all $x,y \in \omega$,
\begin{equation}\label{holder}
\left| \frac{\hat f'(x)}{\hat f'(y)} - 1 \right|
\leq C \beta^{s(x,y)}.
\end{equation}
Moreover the tail \( |\{R>n\}| \) of the return times satisfy the
following estimates:

 \begin{description}
     \item [Summable case]  Under no conditions on $d_n(c)$ 
other than which stem from ($*$)
     \[
     \sum_{n}|\{ R > n\}| < \infty.
     \]
 \item [Polynomial case]
 If $d_n(c) \leq C n^{-\alpha}$
  for all $c \in \Crit$ and $n \geq 1$, then there exists
$\tilde C > 0$ such that
  $$|\{ R > n \}| \leq \tilde C n^{-\alpha}.
  $$
 \item [Stretched exponential case]
 If $b_n(c) \leq C e^{-\beta n^{\alpha}}, \ \alpha\in
(0,1), \
 \beta > 0$
 for all $c \in \Crit$ and $n \geq 1$,
 then for each $\tilde \alpha\in (0,\alpha)$
 there exist $\tilde \beta, \tilde C > 0$ such that
 $$|\{ R > n \}| \leq \tilde C e^{-\tilde \beta n^{\tilde\alpha}}.$$
 \item [Exponential case]
 If $b_n(c) \leq C e^{-\beta n}$, $\beta >0$ for all $c
\in \Crit$
 and $n \geq 1$, then
 there exist $\tilde \beta, \tilde C > 0$ such that
 $$|\{ R > n \}| \leq \tilde C e^{-\tilde \beta n}.$$
 \end{description}

\end{proposition}

In Section~\ref{full returns} we explain how to choose \( \Omega_{0}
\) and how to define the partition \( \mathcal Q \) and the return
time function \( R \). Notice that \( R \) is not a first return time.
In Section~\ref{bounded distortion} we prove the distortion bound and
in Section~\ref{return time estimates} we prove the estimates on the
return times.

\subsection{Large scales and full returns}
\label{full returns}

Let \( \Omega_{0} \subset \Delta \)
be a small neighbourhood of a point $c \in \Crit$ (the precise
requirements on its size will be given in the proof of
Lemma~\ref{fixedprop} below). 
Let \( J\subset X \) be an arbitrary interval.
Consider the map \( F=f^{\hat p}: J \to X \) and the associated
partition \( \hat{\mathcal P}\) on
\( J \) with the stopping time function
\( \hat p \) as defined in Section~\ref{large scales}.

\begin{lemma}\label{fixedprop}
    There exist $t_0 \in {\mathbb N}$ and $\xi>0$
    independent of \( J \)
    such that for every \( \omega\in\hat{\mathcal P} \) there exists
    \( \tilde\omega\subset\omega \) satisfying the following
    properties:
    \begin{itemize}
 \item
  \( f^{\hat p(\omega) + t } \)
 maps \( \tilde\omega \) diffeomorphically onto \( \Omega_0 \)
for some $t \leq t_0$;
 \item
 \( |\tilde\omega|\geq \xi |\omega| \).
\item both components of $f^{\hat p(x)}(\omega \setminus \tilde \omega)$
have length $\geq \delta'/3$.
    \end{itemize}
\end{lemma}

\begin{proof}
By definition of $X$, the preimages of $c$ are dense in $X$.
Therefore there exists \( t_0 \geq 1 \) such that every interval of length
$\geq \delta'$ contains a point $x \in \cup_{t \leq t_0} f^{-t}(c)$
in its middle fifth. Say $f^t(x) = c$.
Now choose
    sufficiently small neighbourhoods \( \omega_x \)
    of each such $x$
    not containing any points of
    \( f^{-j}(\Crit) \) for any \( j < t \). Clearly \( f^t \)
    maps $\omega_x$ diffeomorphically
    to some critical neighbourhood.  By adjusting the size
    of \( \omega_x \) we can make sure
    that they all (i.e. for all points $x$) map onto exactly
    the same critical neighbourhood \( \Omega_{0} \) and
    that $|\omega_x| \leq \delta'/15$.
Let $\tilde \omega \subset \omega$ be the interval that is mapped
onto $\omega_x$ by $f^{\hat p(\omega)}$.
This proves the first and third statement.

 From Lemma~\ref{distortion hat p} we know that the distortion
$f^{\hat p(\omega)}|\omega$ is bounded by $K = K(\epsilon)$.
%
The second statement follows immediately.
\end{proof}

Having fixed $\Omega_0$, let $\delta'' = \min\{ \delta'/3, |\Omega_0|\}$.
In the remainder we will only need to consider intervals $J$ of 
size $\geq \delta''$.

We now define \( \hat
f:\Omega_{0}\to\Omega_{0} \), the associated partition \( \mathcal Q \)
and the stopping time function \( R \) constant on elements of \(
\mathcal Q \) such that \( \hat f= f^{R(\omega)} \) on \(
\omega\in\mathcal Q \).
For each \( \omega \) in the partition \( \hat{\mathcal P}
 \) of \( \Omega_{0} \), let \( \tilde\omega \) denote the
 subinterval given in Lemma~\ref{fixedprop}, so
$|f^{\hat p(\omega)}(\tilde \omega)| \geq \delta'$.
We put  \( \tilde\omega\in\mathcal Q \)
 by definition and \( R(\tilde\omega)= \hat p(\omega)+ t \).
Both components of \( f^{\hat
p(\omega)}(\omega) \setminus \omega_x \)
have size at least $\delta'/3$.
Considering them as new starting intervals
we carry out the construction of Section~\ref{large scales} and repeat the
procedure described above. This
determines all the necessary objects.
In this way each \( \omega\in\mathcal Q \) also has an
associated sequence of \emph{large
scale times} before a
\emph{full return}. We write \( \hat p_{1}=\hat p(x) \)   and
\( \hat p_{i+1}(x) = \hat p_{i} (x) + \hat p (f^{\hat p_{i}(x)}(x))
\) so that \( \hat p_{i+1}(x) \) denotes the total number of iterates
making up the first \( i+1\) large scale stopping times associated to
the point \( x \).  We have \( R(\omega)= \hat
p_{s}(\omega) + t \) for some \( s\geq 1, t \leq t_0 \).

We prove two easy but important consequences of the construction.

\begin{lemma}\label{koebe}
    For each \( n\geq 0 \) and each interval \( \omega \) on which \(
    \hat f^{n} \) is continuous, the distortion of \( \hat f^{n}|\omega \)
    is uniformly bounded (independently of \( \omega \)).
\end{lemma}
\begin{proof}
    The statement follows directly from the construction and
    Lemma~\ref{fixedprop}. Indeed, the third item of
Lemma~\ref{fixedprop} shows that the Koebe space around
$f^{\hat p(\omega)}$ is at least $\delta'/3$.
The additional \( t \) iterates do not significantly affect the distortion.
\end{proof}

\begin{lemma}\label{independent}
For every $i$,
    \[
    |\{x;
   \hat p_{i+1}(x) \text{ exists and }\hat p_{i+1}> \hat p_{i}+ k | \hat
p_{i}
\}|
\leq
 \frac{3K}{\delta'} |\{\hat p > k\}|.
    \]
    Here the expression on the left denotes the conditional
    probability of \( \hat p_{i+1}> \hat p_{i}+ k \) on
    the set of intervals on which \( \hat p_{i} \) is defined.
\end{lemma}

\begin{proof}
    The statement follows immediately from Lemmas~\ref{fixedprop} and
    \ref{koebe}. Indeed,
let $\omega$ be a maximal interval on which $\hat p_i$ is defined and
constant, say $f^{\hat p_i(\omega)}(\omega) = J \supset \omega_x$,
where $\omega_x$ is as in Lemma~\ref{fixedprop}.
Let $\omega' \subset \omega$ be such that $J' = f^{\hat
p_i(\omega)}(\omega')$
is a component of $J \setminus \omega_x$.
By construction $|J'| \geq \delta'/3$.
As the transformation $f^{\hat p_i(\omega)}|\omega$
has distortion bounded by $K = K(\epsilon)$, we get
$$
|\{ x \in \omega'; \hat p_{i+1}(x) > \hat p_i(x)+k \} |
\leq K \frac{|\omega'|}{|J'|}
|\{ y \in J'; \hat p_{J'}(y) > k\}|
$$
Because $|\{ \hat p_{J'} > k \} | \leq | \{ \hat p > k \} |$
the result follows by summing over all the intervals $\omega'$.
\end{proof}

\subsection{Bounded distortion}
\label{bounded distortion}

The function $s$ from \eqref{separation} is called the
\emph{separation time function}.
Notice that \( s(x, y) \) is finite
for all \( x \neq y \), because otherwise $f^n|(x,y)$ would be
homeomorphic for all $n$. The assumptions
on $D_n$ imply that $|(f^n)'(x)|$ does not converge to $0$ for any
$x \in X \setminus \cup_n f^{-n}(\Crit)$, so this cannot happen.
By the same token one can show that some iterate of $\hat f$
is uniformly expanding, i.e. there exists $N$ such that
$|(\hat f^N)'(x)| \geq 2$ wherever it is defined.

\begin{lemma}\label{holder property}
There exists $\beta \in (0,1)$ and $C > 0$
such that for all $\omega \in {\mathcal Q}$ and all $x,y \in \omega$,
\begin{equation}\label{holder2}
\left| \frac{\hat f'(x)}{\hat f'(y)} - 1 \right|
\leq C \beta^{s(x,y)}.
\end{equation}
\end{lemma}

\begin{proof}
For small values of $s(x,y)$, \eqref{holder2} follows immediately
from Lemma~\ref{koebe}
Otherwise, uniform expansion of $\hat f^N$ and Lemma~\ref{koebe}
imply that
$|\hat f(x) - \hat f(y)| \leq |\Omega_0| K 2^{-s(x,y)/N}$.
Because the Koebe space around $\hat f|\omega$ is at least
$\delta'/3$, we get
\begin{eqnarray*}
\left| \frac{\hat f'(x)}{\hat f'(y)} - 1 \right|
&\leq& \left| \left( \frac{\delta'/3|\hat f(x)- \hat f(y)| + 1}
{\delta'/3|\hat f(x)- \hat f(y)|} \right)^2 -1 \right| \\
&\leq&
\left| \left( \frac{\delta' + 3|\Omega_0| K
2^{-s(x,y)/N}}{\delta'} \right)^2 -1 \right|
\leq C2^{-s(x,y)/N},
\end{eqnarray*}
where $C = \frac{6 K |\Omega_0|}{\delta'} + (\frac{3|\Omega_0|
K}{\delta'})^2$.
Here we used $K(\delta'/3) = \left( \frac{1+\delta'/3}{\delta'/3} \right)^2$
as Koebe distortion constant, see \cite[Chapter IV]{MS} in the negative
Schwarzian case. In the general case, we take the constant
from Theorem B in \cite{vSV}.
\end{proof}

\subsection{Return time estimates}
\label{return time estimates}

We fix \( n\geq 1 \) and consider the tail \( \{R>n\} \) of the return times
for \( \hat
f \) on \( \Omega_{0} \). Let us agree to use the notation 
$|\{\hat p > n\}| := 
\sup\{ |\{ x \in J ; \hat p(x) > n\} |/|J| ; |J| \geq \delta''\}$,
which was estimated in Proposition~\ref{tail of hat p}.
In the summable case, no explicit estimates were given, except
that $\sum_n |\{\hat p > n\}| < \infty$.

Before starting the proof we introduce some notation.
Recall that by construction each \( \omega\in\mathcal Q \) has an
    associated sequence
    \[
    0= \hat p_{0}<\hat p_{1} < \hat p_{2} < \dots < \hat p_{s
    (\omega)}< R(\omega)
    \]
    with \( R(\omega) = \hat
    p_{s(\omega)}+t \) and clearly \( s\leq R \).  Write
    \( \mathcal Q^{(n)}=\{\omega\in\mathcal Q ; R(\omega)>n\}\) and let
    \[
    \mathcal Q^{(n)}_{i}=\{\omega\in\mathcal Q^{(n)} ;
\hat p_{i-1}< n \leq \hat p_{i}\}
    \]
    denote the set of elements of \( \mathcal Q \) with \( R(\omega)>n
    \) and having exactly \( i-1 \) large scale
    times before time \( n \).  Moreover, for each \( i \) and every
    sequence \( (k_{1}, \dots, k_{i}) \) of positive integers with
    \( \sum k_{j}=n \) we write
    \[
    \mathcal Q^{(n)}_{i}(k_{1}, \dots, k_{i}) =
    \{\omega\in\mathcal Q^{(n)}_{i} ;
    k_j = \hat p_{j}-\hat p_{j-1} \text{ for }
    \ j\leq i-1, \ k_i = n-\hat p_{i-1}\}.
    \]
    Finally we let
    \[
    |\mathcal Q^{(n)}_{i}| = \sum_{\omega\in\mathcal Q^{(n)}_{i}}|\omega|
    \quad\text{and}\quad
    |\mathcal Q^{(n)}|= \sum_{i\leq n}|\mathcal Q^{(n)}_{i}|.
    \]
    Obviously \( |\{R>n\}|=|\mathcal Q^{(n)}| \).  We are now ready to
    prove Proposition~\ref{return times}.

\begin{proof}[Proof of Proposition~\ref{return times}]
In the stretched exponential case,
take $\tilde \alpha < \hat \alpha < \alpha$, where $\hat \alpha$
is as in Proposition~\ref{tail of hat p}. Both $\tilde \alpha$ and
$\hat \alpha$ can be arbitrarily close to $\alpha$.
In the exponential case take $\tilde \alpha = \alpha = 1$.
    Let \( \eta\in (0,1) \)  be a small number to be determined
    below, depending on \( \alpha \) and \( \beta \) but not on \( n
    \).  We write
    \begin{equation}\label{split}
        |\{R>n\}|=\sum_{i \leq n} |\mathcal Q_{i}|
     = \sum_{i < \eta n^{\tilde \alpha}} |\mathcal Q_{i}|+
\sum_{\eta n^{\tilde \alpha} \leq i \leq n} |\mathcal Q_{i}|.
    \end{equation}
Lemma~\ref{fixedprop} says that a fixed proportion
    \( \xi \)
    of every element in \( \mathcal Q^{(n)}_{i-1} \) has a full return
    to \( \Omega_{0} \) before its next large scale time. Therefore
    \[
    |\mathcal Q^{(n)}_{i}|/|\mathcal Q^{(n)}_{i-1}|\leq 1-\xi.
    \]
    This implies
    \(
    |\mathcal Q^{(n)}_{i}| \leq (1-\xi)^{i}
    \)
    and therefore the second term in \eqref{split}
    satisfies
    \begin{equation}\label{split2}
 \sum_{\eta n^{\tilde \alpha} \leq i \leq  n} |\mathcal Q^{(n)}_{i}|
 \leq  \sum_{\eta n^{\tilde \alpha} \leq i \leq  n}  (1-\xi)^{i} \leq
 \frac1{\xi} (1-\xi)^{\eta n^{\tilde \alpha}}.
    \end{equation}
For the first term write
\[
\sum_{i < \eta n^{\tilde \alpha}} |\mathcal Q^{(n)}_{i}|=
\sum_{i<\eta n^{\tilde \alpha}} \sum_{
\substack{
(k_{1}, \dots, k_{i}) \\ \sum k_{j}=n}}
|\mathcal Q^{(n)}_{i}(k_{1}, \dots, k_{i})|.
\]
For a given sequence \( (k_{1}, \dots, k_{i}) \), Lemma~\ref{independent}
and Proposition~\ref{tail of hat p} imply
\begin{align*}
       |\mathcal Q^{(n)}_{i}(k_{1}, \dots, k_{i})|  &\leq
       |\{\hat p_{i}>\hat p_{i-1}+k_{i-1}-1| \hat p_{i-1}\}|
\cdot \cdot \cdot
       |\{\hat p_{1}>k_{i-1}-1\}| \\
       &\leq
       \tilde K^{i} \prod_{j=1}^{i} |\{\hat p > k_{j}-1\}| \\
       & \leq \tilde K^{i} \prod_{j=1}^{i} e^{-\hat \beta
       (k_{j}-1)^{\hat \alpha}}
\leq (\tilde K e^{\hat \beta})^i e^{-\hat \beta n^{\hat \alpha} }
   \end{align*}
   Here \( \tilde K = 3K/\delta' \) is the constant in the
   statement of Lemma \ref{independent}.
{}From Lemma~\ref{Stirling} we have that the number of
sequences \( (k_{1}, \dots, k_{i}) \) as above equals
$N^+_{n,i}$ and satisfies
$$
N^+_{n,i} \leq
\left\{ \begin{array}{ll}
e^{\hat \eta n^{\tilde \alpha} \log n}  &\text{ if }
     i \leq \eta n^{\tilde \alpha}, \ \tilde \alpha < 1, \\
e^{\hat \eta n} &\text{ if } i \leq \eta n,
\end{array} \right.
$$
for some $\hat \eta = \hat \eta(\eta, \tilde \alpha)$
tending to $0$ as $\eta \to 0$. In the stretched exponential case
\begin{align*}
\sum_{i < \eta n^{\tilde \alpha}} |\mathcal Q_{i}| &=
\sum_{i<\eta n^{\tilde \alpha}} \sum_{
\substack{(k_{1}, \dots, k_{i}) \\ \sum k_{j}=n}}
|\mathcal Q^{(n)}_{i}(k_{1}, \dots, k_{i})|  \\
& \leq \sum_{i< \eta n^{\tilde \alpha}}
e^{\hat \eta n^{\tilde \alpha} \log n} (K e^{\hat \beta})^{i}
e^{-\hat \beta n^{\hat \alpha}}
 \leq \hat C e^{-\beta' n^{\hat \alpha}}
\end{align*}
for some \( \hat C, \beta' >0 \) as long as \( \hat \eta  \) is
sufficiently small. In precisely the same way we get
$\sum_{i < \eta n} |\mathcal Q_i| \leq \hat C e^{-\beta' n}$
in the exponential case.

 To treat the summable and polynomial case we write \( |\{R>n\}| \) as in
 \eqref {split}, with \( \hat \alpha =1 \) and \( \eta = 1/2 \).
The same argument gives an exponential estimate as in
 \eqref {split2} for the second term.  To estimate the
first term, notice that for each \( i \) and each
sequence \( k_{1}, \dots, k_{i}
 \) with \( \sum k_{j}=n \), the largest $k_j$ satisfies
 \( k_{j}\geq n/i \).  Thus letting
 \[
 \mathcal Q^{(n)}_{i, j}= \{\omega\in\mathcal Q^{(n)}_i ;
 k_{j'}< n/i \text{ for } j'<j \text{ and } k_{j}\geq n/i\}
 \]
 we have
\begin{eqnarray*}
\sum_{n \geq 1} \sum_{i < n/2} |\mathcal Q^{(n)}_{i}| 
&=&
\sum_{n\geq 1} 
\sum_{i <  n/2} \sum_{j=1}^{i} |\mathcal Q^{(n)}_{i, j}| \\
&\leq& \sum_{n \geq 1} \sum_{i < n/2} i (1-\xi)^{i-1} |\{ \hat p > n/i \}| \\
&\leq& \sum_{i \geq 1} i (1-\xi)^{i-1} 
\sum_{n \geq 1} |\{ \hat p > n/i \}|.
\end{eqnarray*}
Substituting $k = \lfloor n/i \rfloor$, and using the fact that
at most $i$ different values of $n$ give the same value of $k$, we find
that the above is bounded by
$\sum_{i \geq 1} i^2 (1-\xi)^{i-1} \sum_{k \geq 1}  |\{ \hat p > k \}|$
which is finite (use Proposition~\ref{tail of hat p}).

In the polynomial case we get
 \begin{equation*}
  \sum_{i < n/2} |\mathcal Q^{(n)}_{i}|
 \leq
{\mathcal O}(n^{-\alpha}) \sum_{i < n/2} i^{1+\alpha} (1-\xi)^i
= {\mathcal O}(n^{-\alpha}).
\end{equation*}
Together with the exponential estimate for the
term $\sum_{i \geq n/2} |{\mathcal Q}^{(n)}_i|$, this yields the
proposition.
\end{proof}

\subsection{Proof of Theorems~\ref{existence}, \ref{decay} and \ref{CLT}}
\label{final_proofs}
We now state the assumptions and results of Young which we want to
apply. Together with the estimates obtained in
Proposition~\ref{return times}, they easily
imply Theorems~\ref{existence}, \ref{decay} and \ref{CLT}.
Let $m$ denote Lebesgue measure on $X$.
L.-S. Young applies the following tower construction for her results.
Given a countably piecewise monotone and onto map
$\hat f: \cup_{\omega \in {\mathcal Q}} \omega \to \Omega_0$,
$\hat f|\omega = f^{R(\omega)}$,
define a tower
$$
\Omega =
\bigsqcup_{ \substack{ \omega \in {\mathcal Q} \\ 0 \leq i < R(\omega)}}
 (\omega,i),
$$
with an action
$$
g(x,i) = \left\{ \begin{array}{ll}
(x,i+1) & \text{ if } x \in \omega, \, i+1 < R(\omega), \\
(\hat f(x), 0) & \text{ if } x \in \omega, \, i+1 = R(\omega).
\end{array} \right.
$$
The connection with the original map $f$ is established by means of
the projection $\pi(x,i) = f^i(x)$.
Because $f^i$ is smooth and has bounded distortion on each
$\omega \in {\mathcal Q}$, $i < R(\omega)$,
this projection has bounded distortion.
Also $\pi \circ g = f \circ \pi$.
Therefore, if $\nu$ is a $g$-invariant absolutely continuous
probability measure
on $\Omega$, $\mu := \nu \circ \pi^{-1}$ is an invariant
absolutely continuous probability
measure on the interval.

We summarize Young's results from \cite{Y} as far as we need them.
For a fixed $\beta \in (0,1)$ as in Lemma~\ref{holder property},
let
$$
{\mathcal C_{\beta}} = \{ \phi:\Omega \to {\mathbb R};
\, \exists \,  C > 0 \, \forall \, x,y \,\,
|\phi(x)-\phi(y)| \leq C\beta^{s(x,y)} \}
$$
and
$$
{\mathcal C}_{\beta}^+ = \{ \phi \in {\mathcal C_{\beta}} ;
\phi \geq 0\}.
$$
Here we have extended the separation time $s$ to $\Omega$ in the obvious
way.
Also let $m_{\Omega}$ be Lebesgue measure on $\Omega$.
(A priori, $m_{\Omega}$ can be infinite.)

\begin{theorem*}[Young  \cite{Y}]\label{young}
Suppose that $\hat f:\Omega_0 \to \Omega_0$ is as above,
i.e.
$m(\Omega_0 \setminus \cup_{\omega \in {\mathcal Q}} \omega) = 0$
and \eqref{holder2} holds.
Let $\{ \rho_n \}$ be a sequence of positive reals related to the
tail behaviour of $R$ as follows.
If $m(\{ R > n \}) \leq n^{-\alpha}$, then $\rho_n = n^{1-\alpha}$,
if $m(\{ R > n \}) \leq e^{-\beta n}$, then $\rho_n = e^{-\beta' n}$
for some (any) $\beta' < \beta$ and
if $m(\{ R > n \}) \leq e^{-n^{\alpha}}$ for some $\alpha \in (0,1)$,
then $\rho_n = e^{-n^{\alpha'}}$ for some (any) $\alpha' < \alpha$.
Then
\begin{enumerate}
\item
If $\sum_n m(\{ R > n\}) < \infty$, then
$\Omega$ carries an $g$-invariant
absolutely continuous probability measure $\nu$ (Kac's Theorem)
and
$\frac{d\nu}{dm_{\Omega}} \in {\mathcal C}_{\beta}^+$.

\item
For any measure $\tilde \nu$ with
$\frac{d\tilde \nu}{dm_{\Omega}} \in {\mathcal C}_{\beta}^+$,
$g_*^n \tilde \nu \to \nu$ and there exists $C_{\tilde \nu} > 0$
such that
$|g_*^n \tilde \nu - \nu| \leq C_{\tilde \nu} \rho_n$.

\item
For any pair of functions
$\phi \in L^{\infty}(\Omega,m_{\Omega})$ and
$\psi \in {\mathcal C}_{\beta}$, there exists $C_{\phi,\psi} > 0$ such that
$$
|\int (\phi \circ g^n) \psi d\nu - \int \phi d\nu \int \psi d \nu|
\leq C_{\phi,\psi} \rho_n.
$$
\item
If $m(\{ R > n\}) \leq {\mathcal O}(n^{-\alpha})$ for some $\alpha > 2$,
then
for any $\phi \in {\mathcal C}_{\beta}$ which is not a coboundary
($\phi \neq \psi \circ g - \psi$ for any $\psi$),
the Central Limit Theorem holds, i.e.
there exists $\sigma > 0$ such that
$\frac1{\sqrt n} \sum_{i=0}^{n-1} \phi \circ g^i$ converges to
the normal distribution ${\mathcal N}(\int \phi d\nu, \sigma)$.
\end{enumerate}

\end{theorem*}

{\em Remark: }
Young states this theorem in terms of a stopping time $\hat R$
which is the extension of $R$ to the entire tower $\Omega$.
As it happens $m_{\Omega}(\{ \hat R > n \}) =
\sum_{k \geq n} m(\{ R > k \})$, so that
$m_{\Omega}(\{ \hat R > n \}) \leq {\mathcal O}(n^{-\alpha})$
if $m(\{ R > k \}) \leq {\mathcal O}(n^{-\alpha-1})$.
This explains why the exponent in the polynomial
case at first glance looks different from the ones in Young's
version.

\medskip

Using the projection $\pi$, these results immediately carry over to
the original map $f$ with measure $\mu = \nu \circ \pi^{-1}$.
   Using the projection $\pi$, immediately carry over to
the original map $f$ with measure $\mu = \nu \circ \pi^{-1}$.
With respect to the support of the measure, note that
$\hat f: \Omega_0 \to \Omega_0$
is a mixing map, and its invariant measure
$\frac1{\nu(\Omega_0)} \nu|\Omega_0$ has the whole interval $\Omega_0$
as support. The formula $\mu = \nu \circ \pi^{-1}$ shows
that $\Omega_0 \subset \mbox{supp}(\mu)$.

Finally, recall from Proposition~\ref{return times} how the tail
$m(\{ R > n \})$ is related to $d_n(c)$. Therefore Theorem~\ref{young}
immediately gives Theorems~\ref{existence}, \ref{decay} and \ref{CLT}.

\medskip
\noindent
Department of Mathematics\\
University of Groningen\\
P.O. Box 800, 9700 AV Groningen\\
The Netherlands\\
\texttt{bruin@math.rug.nl}\\
\texttt{http://www.math.rug.nl/\~{}bruin}

\medskip
\noindent
Department of Mathematics\\Imperial College\\
180 Queen's Gate, London SW7\\
UK\\
\texttt{Stefano.Luzzatto@ic.ac.uk} \\
\texttt{http://geometry.ma.ic.ac.uk/\~{}luzzatto}

\medskip

\noindent
Department of Mathematics\\
University of Warwick\\
Coventry CV4 7AL\\
UK\\
\texttt{strien@maths.warwick.ac.uk}\\
\texttt{http://www.maths.warwick.ac.uk/\~{}strien}

\end{document}